\documentclass[a4paper,fleqn]{cas-dc}

\usepackage[numbers]{natbib}

\usepackage[per-mode=symbol]{siunitx}
\usepackage{cases}
\usepackage{graphicx}
\usepackage{tabularx}
\usepackage{caption}
\usepackage{threeparttable}

\def\tsc#1{\csdef{#1}{\textsc{\lowercase{#1}}\xspace}}
\tsc{WGM}
\tsc{QE}

\begin{document}
\let\WriteBookmarks\relax
\def\floatpagepagefraction{1}
\def\textpagefraction{.001}

\shorttitle{Impact of Data Quality on Renewable Energy Potential Estimations}    

\shortauthors{S. Risch, R. Maier, J. Du, N. Pflugradt, L. Kotzur, D. Stolten}  

\title [mode = title]{Impact of Data Quality on Renewable Energy Potential Estimations}  



%

\author[a]{Stanley Risch}[orcid=0000-0002-7188-155X]

\cormark[1]


\ead{s.risch@fz-juelich.de}




\address[a]{Institute of Energy and Climate Research - Techno-economic Systems Analysis (IEK-3), Forschungszentrum Jülich GmbH, Wilhelm-Johnen-Straße, 52428 Jülich, Germany}
\address[b]{RWTH Aachen University, E.ON Energy Research Center, Institute for Energy Efficient Buildings and Indoor Climate, Mathieustraße 10, Aachen, 52074, Germany}
\address[c]{TH Köln, Cologne Institute for Renewable
Energy (CIRE), Betzdorfer Straße 2, Köln, 50679, Germany}
\address[d]{RWTH Aachen University, c/o Institute of Energy and Climate Research - Techno-economic Systems Analysis (IEK-3), Forschungszentrum Jülich
GmbH, Wilhelm-Johnen-Straße, Jülich, 52428, Germany}
\author[a]{Rachel Maier}[]
\author[b]{Junsong Du}[]
\author[a]{Noah Pflugradt}[]
\author[c]{Peter Stenzel}[]
\author[a]{Leander Kotzur}[]
\author[a,d]{Detlef Stolten}[]




\credit{}


\cortext[1]{Corresponding author}
\fntext[1]{Tool for Renewable Energy Potentials - Database\\(doi: \href{https://doi.org/10.5281/zenodo.6414018}{10.5281/zenodo.6414018})}


\begin{abstract}
Potential analyses identify possible locations for renewable energy installations, such as as wind turbines and photovoltaic arrays. 
The results of previous potential studies, however, are not consistent due to different assumptions, methods, and datasets.
In this study, we compare commonly used land use data sources with regard to area and position. Using Corine Land Cover leads to an overestimation of the potential areas in a typical wind potential analysis by a factor of 4.6 and 5.2 in comparison to Basis-DLM and Open Street Map, respectively. 
Furthermore, we develop scenarios for onshore wind, offshore wind, and open-field photovoltaic potential estimations based on land eligibility analyses and calculate rooftop photovoltaic potential using 3D building data. 
The potential capacities and possible locations are published for all administrative levels in Germany in the freely accessible database \textbf{trep-db},\textsuperscript{1} for example, to be incorporated into energy system models.
The investigations are validated using high-resolution regional potential analyses and benchmarked against other studies in the literature. 
Findings from the literature, which can be used by legislators to design regulation, are rarely comparable and consistent due to differences in the datasets used. 

\end{abstract}



\bibliographystyle{model1-num-names}
\ExplSyntaxOn
\keys_set:nn { stm / mktitle } { nologo }
\ExplSyntaxOff
\maketitle



\section{Introduction}\label{sec:intro}
In 2015, the Paris Climate Agreement~\citep{united_nations_paris_2015} was signed by 195 countries with the aim of limiting global temperature increases to below \SI{2}{\celsius} above pre-industrial levels. 
Germany strengthened its ambitions to combat climate change in the Climate Protection Act of 2021~\citep{noauthor_bundes-klimaschutzgesetz_2021} by setting the goal of climate-neutrality in 2045.
To achieve this target, the capacity of renewable energy technologies must be greatly increased \citep{stolten_strategien_2021}. 
However, the specific land requirement of renewable energy technologies as photovoltaic (PV) systems or wind turbines exceeds those of conventional power plants, and therefore larger areas will be needed for the future energy supply. 
The eligibility and availability of construction areas constitute a limiting factor and determine the region-specific potentials of renewable energy sources.
\par
Energy system models and corresponding studies can support the planning process of future energy systems with high shares of renewable energy. 
However, the set potentials or renewable expansion goals significantly differ in national energy system studies~(see \citep{stolten_strategien_2021,luderer_deutschland_2021,brandes_wege_2021,kendziorski_100_2021,prognos_klimaneutrales_2021,fette_multi-sektor-kopplung_2020,goke_100_2021}) and often information about methods, data sources, or assumptions is lacking. 
As energy system analyses can have different regional scopes and levels, internally-consistent potential data on different geographical levels is needed. 

Regionalized potentials are estimated with potential analyses, which are performed sequentially:
First, eligible areas for the construction of renewable power installations are determined. Then on the basis of these areas, the potential capacity and energy generation can be estimated. 
The definition of eligible areas therefore plays a crucial role in the estimation of potentials.
Depending on the technology, the calculation of eligible areas can be performed by means of statistical formulas, with building models, or with land eligibility analyses. 
For technologies requiring 'open space', such as open-field PV and wind, the latter methodology is utilized, which requires geospatial datasets.
The results of the eligible areas and the potential capacity are therefore influenced by these datasets; however, the impact has not yet been evaluated~\citep{mckenna_high-resolution_2021}. 
The comparability of studies using different datasets is therefore unknown.

The following analyses the impact of commonly used datasets on land eligibility analyses and presents scenarios for onshore, offshore wind, open-field PV, and rooftop PV potentials in Germany. 
The structure of the paper is twofold: 
In the first part (\autoref{sec:LEAPotentials}), the potential analyses on open spaces for open-field PV and offshore and onshore wind are presented.
These are based on land eligibility analyses (\autoref{ssec:LEA}) using geospatial datasets.
We first evaluate land use datasets in the context of potential analyses in \autoref{sssec:data}, which closes a gap in the literature and also underlines the quality of the following potential analyses.
Then, we describe the state of the art and methodology before presenting the results and discuss these for the following technologies: onshore wind (\autoref{ssec:onwind}), offshore wind (\autoref{ssec:offwind}), and open-field PV (\autoref{ssec:ofpv}).
In the second part of the paper, we present the potential analysis for rooftop PV in \autoref{sec:rtpv}.
Thereby, for the first time, 3D building data is used to estimate the potential for Germany.
For the considered technologies, we benchmark the results against other potential studies to provide insights into the impact of data sources, exclusion criteria and methodologies. 
The results of our potential analyses are published in the open database \textbf{trep-db}, with utility, for instance, for energy system modelers. 
The area and capacity potentials of high quality and resolution are provided for several scenarios per technology, and for the administrative levels spanning the municipality to national levels in Germany.

\section{Renewable energy potentials on open spaces}\label{sec:LEAPotentials}
The following chapter presents the potential analyses for renewable energy technologies that require open space. 
First, the general methodology and evaluation of the input data for land eligibility analyses is presented in \autoref{ssec:LEA} and \ref{sssec:data}, respectively. 
Afterwards, potential analyses are performed for onshore wind~(\autoref{ssec:onwind}), offshore wind~(\autoref{ssec:offwind}), and open-field PV~(\autoref{ssec:ofpv}).

\subsection{Land eligibility analysis}\label{ssec:LEA}
The construction of renewable power generation sites requires eligible land. 
Therefore, land eligibility analyses based on geospatial data are performed as the first step of open space potential analyses.
Geospatial data is used to consider areas eligible, e.g., bare land, or to extract information to identify ineligible areas due to criteria such as physical constraints, such as steep slopes, or regulatory constraints, such as national parks.
Aside from considering areas ineligible, a setback around these areas may also not be usable depending on the exclusion criteria. 
This can be achieved by adding a buffer around identified land use categories.
For instance, the construction of wind turbines is not only ineligible for land use categories such as settlements or streets but also within a certain distance around these areas. 
To perform the land eligibility analysis, the open source tool \textit{GLAES}~\cite{ryberg_evaluating_2018} was applied, which performs the necessary operations using geospatial datasets to determine eligible areas for renewable energy sites. 
We chose a resolution of \SI{10}{\metre}$\times$\SI{10}{\metre} in the land eligibility analyses to ensure the representation of detailed features for the high-resolution potentials.
However, as the chosen datasets have a significant impact on the results of the land eligibility analysis, they are analyzed in the following section.

\subsection{Evaluation of land use datasets}\label{sssec:data}

To identify geographical features, land eligibility analyses require geospatial datasets that describe land characteristics, e.g., land use, elevation, and water depths.
The quality and resolution of the feature representation in the dataset therefore has an impact on the land eligibility analysis.
For instance, missing represented urban areas in a settlement dataset can underestimate excluded areas within the potential analyses or a coarse resolution can either under- or overestimate areas depending on its classification. 
However, even though the datasets influence potential analyses, \citet{mckenna_high-resolution_2021} state that the impact of employing different data sources has not yet been evaluated. 

In the present section, the first comparison and evaluation of the following four land use datasets for potential analyses is performed for the regional scope of Germany:
\begin{itemize}
    \item Basis-DLM~\citep{geobasisdaten__geobasis-de__bkg_2021_digitales_2021}: Official German dataset with a high positional accuracy (between $\pm$\SI{3}{} and $\pm$\SI{15}{\metre} depending on the feature)
    \item Corine Land Cover (CLC) \citep{copernicus_programme_corine_2018}: Land cover raster dataset with \SI{100}{\metre}$\times$\SI{100}{\metre} resolution
    \item Open Street Map (OSM) \citep{openstreetmap_contributors_open_2021}: User based land cover vector dataset
    \item World Database on Protected Areas (WDPA) \citep{unep-wcmc_iucn_world_2016}: Vector dataset with information on protected areas 
\end{itemize}
According to \citet{mckenna_high-resolution_2021}, the datasets CLC, OSM, and WDPA are commonly used global and continental ones for potential analyses.

Two parameters are used to compare the covered area using dataset $A(dataset)$ for different features.
The first parameter is the Normalized Total Area (NTA) (\autoref{eq:NTA}): It describes the covered area by a dataset compared to the maximum area covered by any other dataset capturing the feature:
    \begin{equation}\label{eq:NTA}
        NTA(dataset_i)=\frac{A(dataset_i)}{\max\limits_{j}(A(dataset_j)}
    \end{equation}
The second parameter is the Intersection over Union (IoU) of two datasets (\autoref{eq:IoU}): 
It describes the area identified by both datasets divided by the area identified when combining the two. 
An Intersection over Union of $100\%$ describes that identified areas of two datasets as equal: 
    \begin{equation}\label{eq:IoU}
    IoU(dataset_{(i,j)})=\frac{A(dataset_i) \, \cap\, A(dataset_j)}{A(dataset_i)\, \cup \, A(dataset_j)}
    \end{equation}
    
The calculation of the parameters is performed with additional standard setback distances for wind potential analyses (see the supplementary material).
Furthermore, a rasterized representation of the datasets with the corresponding setback distance and a resolution of \SI{10}{\metre}$\times$\SI{10}{\metre} is used.

\begin{figure*}
\centering
\includegraphics[width=0.99\textwidth]{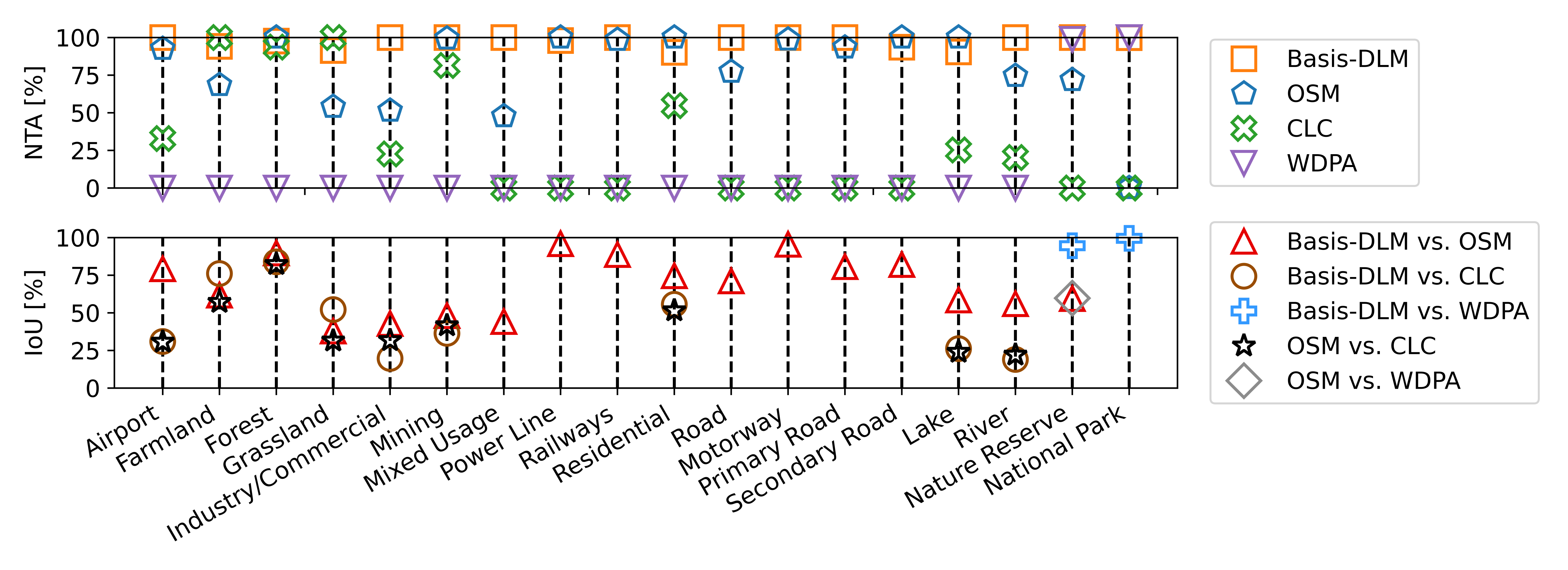}
\caption{Comparison of the Normalized Total Area (NTA) (upper plot) and Intersection over Union (IoU) (lower plot) per category for the covered area with additional setback distances of a wind potential analysis for different data sources.}\label{fig:data_analysis}
\end{figure*}

\autoref{fig:data_analysis} shows the results for the Normalized Total Area and the Intersection over Union for the analyzed land use categories. 
The definition of the filter used to identify the respective features in each dataset can be found in the supplementary information.
\begin{figure}[h]
	\centering
\includegraphics[width=0.49\textwidth]{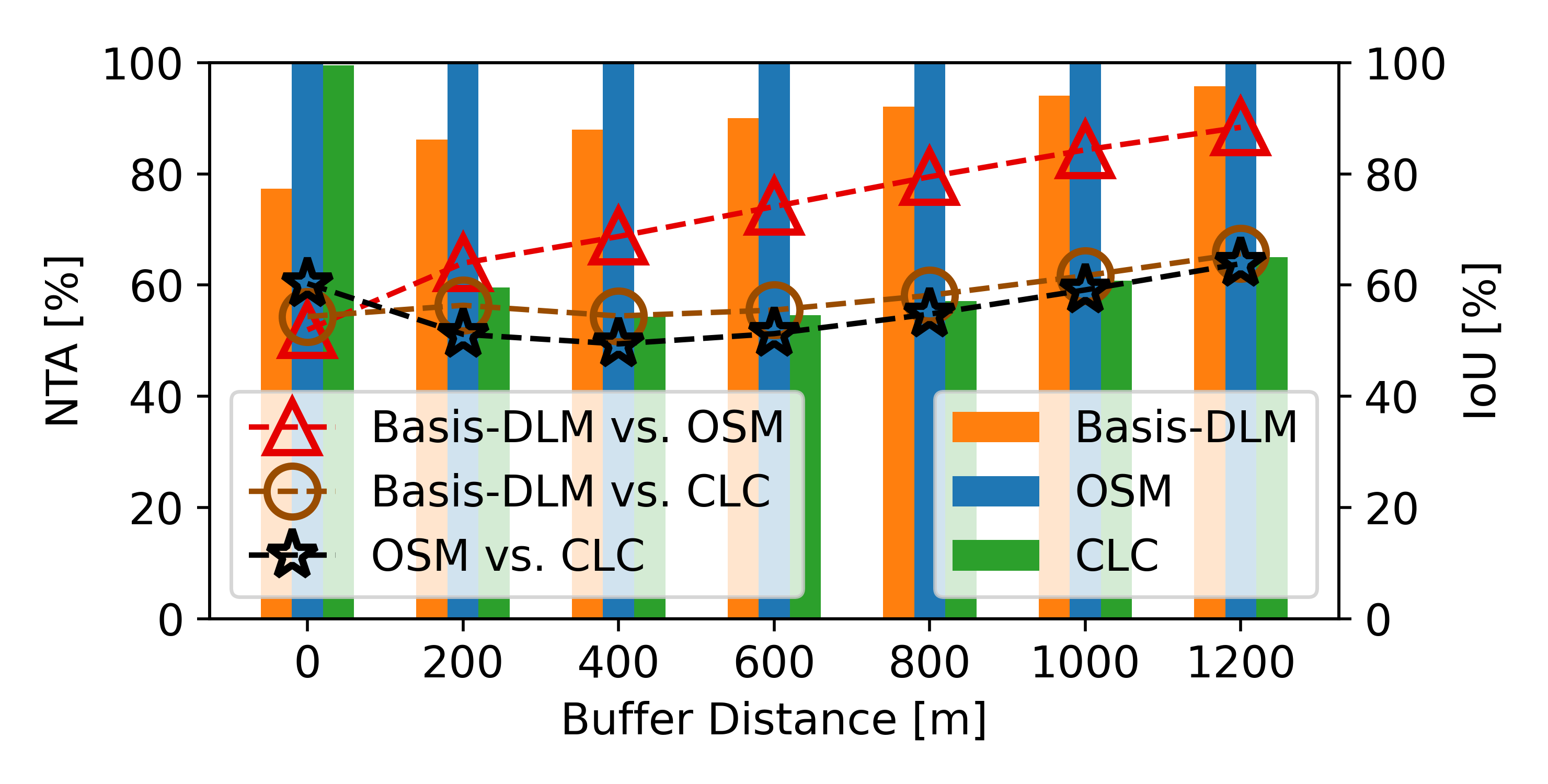}
	  \caption{Intersection over Union (bars) and Normalized Total Area (markers) for the residential category for a variation of the setback distance.}\label{fig:data_analysis_residential}
\end{figure}
It must be stated that the WDPA only provides information on protected areas, in which its results are only comparable for these. 
Furthermore, the CLC does not include all presented categories, such as features with narrow widths like roads or power lines.

For the typical wind potential analysis, OSM and Basis-DLM identify the largest area and share the highest Intersection over Union for most categories.
For categories of larger sizes, i.e., farmland, forests, mining, and grassland, the CLC classifies a similar buffered area as Basis-DLM. Nevertheless, the Intersection over Unions are only high (above $70\%$) for the categories of farmland and forests, which means that the identified areas vary largely in position for the categories of mining and grassland.
In the nature protection categories ('Nature Reserve', 'National Park'),  WDPA and Basis-DLM identify almost identical areas, recognizable through Intersection over Unions near $100\%$. The OSM-identified nature reserves, on the other hand, show significant deviations from the Basis-DLM and WDPA, whereas national parks cannot be identified with the OSM.
\par
\citet{mckenna_high-resolution_2021} and \citet{masurowski_spatially_2016} state that the land use category settlements (often synonymously used with residential land use in land eligibility analyses) highly influences the potential analysis in Germany, which is why the category is further analyzed. 
For a typical setback distance of \SI{1000}{\metre} in \autoref{fig:data_analysis}, the OSM covers the largest area, followed by the Basis-DLM.
\autoref{fig:data_analysis_residential} shows how the Normalized Total Area and Intersection over Union of the residential category behaves when varying the setback distance. Without any buffer, the OSM covers the largest area, closely followed by the CLC. 
With increasing setback distances, the Normalized Total Area of the CLC first drops rapidly, whereas the Normalized Total Areas of the Basis-DLM and OSM converge. 
This indicates that the Basis-DLM and OSM classify more, but smaller settlements in comparison to the CLC. With increasing setback distances of above \SI{600}{\metre}, the Normalized Total Area of the CLC increases again, as the buffered areas around smaller features identified by the Basis-DLM and OSM exhibit increasing overlaps.
The similarity of the Basis-DLM and OSM is also shown by their Intersection over Union, which increases over the buffer distance of up to \SI{88}{\percent}.
The described behavior of CLC in the residential category, on the other hand, indicates that the representation of smaller detached settlements is missing, whereas larger settlements are covered to a coarser extent. 
This can be explained by the relatively coarse resolution of \SI{100}{\metre}$\times$\SI{100}{\metre} of the dataset and the minimum mapping unit of \SI{25}{\hectare}.

Using CLC to classify residential areas can lead to an underestimation of excluded areas and therefore an overestimation of potential areas, especially for typical setback distances for wind potential analyses of between \SI{600}{} and \SI{1000}{\metre}.
Furthermore, CLC cannot represent categories with small size features or widths due to its minimal mapping unit.
Therefore, the usage of CLC for potential analyses is not recommended.
\par
The similarity of OSM and Basis-DLM can be seen for multiple categories, e.g., power lines, motorways, in which the Normalized Total Areas and Intersection over Unions are near $100\%$. 
For these categories, the use of both datasets is justified. 

\begin{figure}[h]
	\centering
\includegraphics[width=0.49\textwidth]{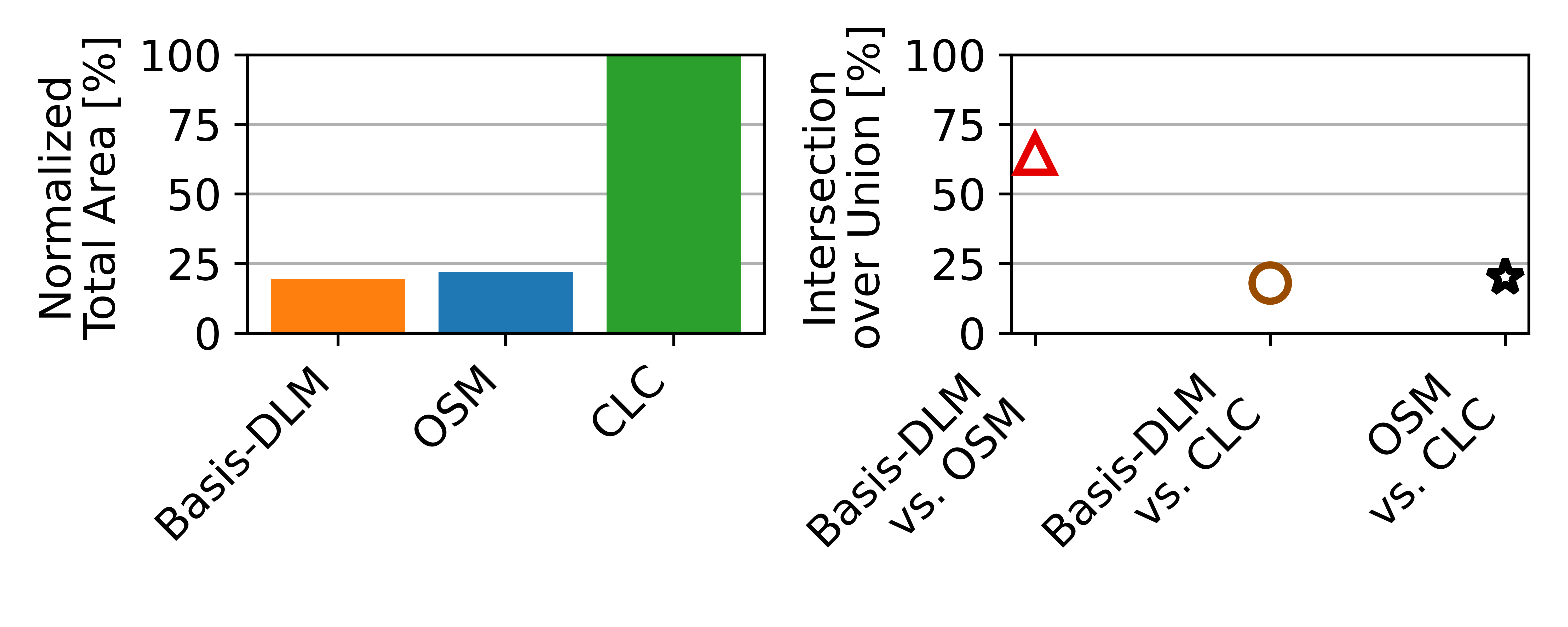}
	  \caption{Normalized Total Area and Intersection over Union for the resulting potential areas of a typical wind analysis using different datasets.}\label{fig:data_analysis_results}
\end{figure}

\autoref{fig:data_analysis_results} shows the results of a land eligibility analysis with the setback distances of a typical wind potential analysis (see the supplementary material) for Basis-DLM, OSM, and CLC. 
When all three of the datasets include the feature, the corresponding dataset is used; otherwise, a default dataset is applied for the exclusion.
In particular, the use of the CLC for land eligibility analyses leads to significantly different results. 
The resulting total potential area using the CLC is roughly $80\%$ higher compared to the OSM and Basis-DLM. 
However, the analyses using Basis-DLM and OSM, even though resulting in similar total areas, have discrepancies locating the areas, which is indicated by an Intersection over Union of $64\%$.

In summary, the presented analysis highlights the discrepancy of land use declaration by different datasets, by which the results of potential analyses using different datasets are hardly comparable. 
Furthermore, the analysis reveals the large sensitivity of land eligibility results regarding the chosen datasets and the importance of data with a high resolution and high positional accuracy for potential analyses. 
\par
For the scenarios outlined in this paper and the enclosed database \textbf{trep-db}, WDPA is used for the nature protection categories and Basis-DLM for most land use categories in the land eligibility analyses. 
Further datasets are employed for other individual categories (e.g., water protection areas), but are not considered in this data analysis due to the lack of comparability to other datasets.
Basis-DLM was selected due to the exhibited characteristics and is preferred over OSM due to its official character, high positional accuracy, detailed documentation, and the categories covered.
Furthermore, for the wind potential analysis, the residential areas are split into inner and outer areas (see \autoref{sssec:onwind-met}), which is not possible by the sole use of OSM.

\subsection{Onshore wind potential}\label{ssec:onwind}

\subsubsection{Literature}\label{sssec:onwind-lit}
Onshore wind analyses are generally performed with a greenfield approach (see, e.g., \citep{ruiz_enspreso_2019,trondle_home-made_2019,ebner_regionalized_2019,peters_raumlich_2015,landesamt_fur_natur_umwelt_und_verbraucherschutz_nordrhein-westfalen_lanuv_potenzialstudie_2021,amme_photovoltaik-_2021,ryberg_uniformly_2020,lutkehus_potenzial_2013,landesanstalt_fur_umwelt_baden-wurttemberg_potenzialanalyse_2019,wiehe_nothing_2021,luderer_deutschland_2021}). 
Only \citet{masurowski_spatially_2016} proceed differently: First, initially available areas (pre-selected areas), among others grassland and arable land, are selected and then further exclusion criteria are applied. 
For the considered studies, apart from the varying approach, the most influential differences correspond to the exclusion definitions, i.e., the used criteria and setback distances, the used data sources (cf. \autoref{sssec:data}), and the capacity estimations. \par
All of the cited studies pertain to residential areas. Nevertheless, the setback distances and the data sources vary significantly.
\citet{trondle_home-made_2019} use the European Settlement Map \cite{european_commission_joint_research_centre_european_2016} to identify built-up areas and exclude them without a buffer, whereas \citet{peters_raumlich_2015} and LANUV~\citep{landesamt_fur_natur_umwelt_und_verbraucherschutz_nordrhein-westfalen_lanuv_potenzialstudie_2021} differentiate the setback distance between inner (\SI{1000}{\metre}) and outer areas (\SI{750}{}, \SI{720}{\metre}) for their analysis based on Basis-DLM, and \citet{amme_photovoltaik-_2021} use OSM~\cite{openstreetmap_contributors_open_2021} to identify residential areas and apply varying setback distances (\SI{400}{} to \SI{1000}{\metre}). 
\citet{masurowski_spatially_2016} specify \SI{1000}{\metre} for housing areas but vary the setback distance in their scenarios, without stating a data source.
Others use CLC~\cite{copernicus_programme_corine_2018} to exclude settlements with \SI{800}{\metre} \cite{ryberg_uniformly_2020,mckenna_cost-potential_2014} to \SI{1000}{\metre} \cite{sensfus_langfristszenarien_2021-2}. 
\citet{ryberg_uniformly_2020} additionally exclude urban areas from EuroStat Urban~\citep{eurostat_eurostat_2011} with \SI{1200}{\metre} setback.
\citet{lutkehus_potenzial_2013} exclude residential areas with \SI{600}{\metre} based on DLM250~\citep{geobasisdaten__geobasis-de__bkg_2012_digitales_2012}.
Moreover, area-based residential, settlement, and urban exclusions of individual residential buildings are specified in two regional potential analyses \cite{landesamt_fur_natur_umwelt_und_verbraucherschutz_nordrhein-westfalen_lanuv_potenzialstudie_2021,landesanstalt_fur_umwelt_baden-wurttemberg_potenzialanalyse_2019} and \cite{luderer_deutschland_2021}.
\par
Another frequently discussed topic is the construction of wind turbines in forests. 
The missing consensus in legislation of the federal states can also be observed in the literature, i.e., Thuringia generally forbids the use of forests for wind turbines, whereas in other federal states, construction of these in forests is ongoing \cite{fachagentur_wind_an_land_entwicklung_2021}. \citet{ruiz_enspreso_2019} exclude forests in the land eligibility analysis, whereas others \cite{trondle_home-made_2019,ryberg_uniformly_2020} allow for construction in forests.
Other studies treat forests in a differentiated manner.
Coniferous forests in densely-wooded municipalities are cited in the LANUV report \citep{landesamt_fur_natur_umwelt_und_verbraucherschutz_nordrhein-westfalen_lanuv_potenzialstudie_2021}. 
\citet{amme_photovoltaik-_2021} and \citet{peters_raumlich_2015} consider scenarios that include and exclude forests. 
Meanwhile, \citet{ebner_regionalized_2019} exclude forests in one scenario and allow the construction in $10\%$ for another.
\citet{wiehe_nothing_2021} successively reduce the usable parts of forests in their more restrictive scenarios. 
\citet{lutkehus_potenzial_2013} exclude forests in federal states with less than $15\%$ forest shares. Furthermore, selective exclusions are performed based on the function of the forest.
\par
Inconsistencies can also be seen for protected landscapes, which are a highly influential exclusion ($28\%$ of Germany's area \cite{unep-wcmc_iucn_world_2016}).
Several studies \cite{peters_raumlich_2015,amme_photovoltaik-_2021,wiehe_nothing_2021,lutkehus_potenzial_2013} regard protected landscapes differently between their scenarios, whereas others
\cite{landesanstalt_fur_umwelt_baden-wurttemberg_potenzialanalyse_2019,ryberg_uniformly_2020, mckenna_cost-potential_2014} generally rule out the construction of wind turbines in protected landscapes. 
\par
After determining the eligible areas, a few studies further reduce it, either by a certain share \citep{trondle_home-made_2019,ebner_regionalized_2019} or based on suitability factors~\citep{mckenna_cost-potential_2014,sensfus_langfristszenarien_2021-2} per land use category of CLC~\citep{copernicus_programme_corine_2018}.
\par
Based on the identified eligible areas, the capacity can be estimated. 
To this end, two different methods are utilized:
The majority of the analyzed studies \cite{landesamt_fur_natur_umwelt_und_verbraucherschutz_nordrhein-westfalen_lanuv_potenzialstudie_2021, landesanstalt_fur_umwelt_baden-wurttemberg_potenzialanalyse_2019, ebner_regionalized_2019, wiehe_nothing_2021, luderer_deutschland_2021,masurowski_spatially_2016,lutkehus_potenzial_2013} 
distribute turbines to eligible areas using a spacing distance specified by a multiple of the rotor diameter (D) in prevailing and transverse wind directions between the individual turbines. 
The turbine spacing ranges from \SI{5}{D} to \SI{9}{D} in the prevailing wind direction and \SI{3}{D} to \SI{4}{D} in the transverse direction.
In contrast, others use fixed capacity densities, e.g., \SI{5}{\mega\watt\per\kilo\metre\squared}~\cite{ruiz_enspreso_2019}, \SI{8}{\mega\watt\per\kilo\metre\squared}~\cite{trondle_home-made_2019}, and \SI{21}{\mega\watt\per\kilo\metre\squared}~\cite{amme_photovoltaik-_2021}.
\citet{sensfus_langfristszenarien_2021-2} and \citet{mckenna_cost-potential_2014} vary the capacity density regionally.

\subsubsection{Methodology}\label{sssec:onwind-met}
The available area for onshore wind was determined with a greenfield land eligibility analysis (see \autoref{ssec:LEA}). To this end, five different scenarios are defined:
\begin{itemize}
    \item S1 Legislation: The exclusions are defined according to the laws of Germany's federal states based on \cite{fachagentur_wind_an_land_uberblick_2021} and own corrections
    \item S2 Expansive: Wind expansion favoring exclusions including forests and protected landscapes\\
    S2a No Protected Landscapes: S2, excluding protected landscapes \\
    S2b No Forests: S2, excluding forests 
    \item S3 Restrictive: Restrictive exclusions
\end{itemize}
In scenario~1, federal state-specific exclusions are applied. The exclusions for each federal state are defined in accordance with \citet{fachagentur_wind_an_land_uberblick_2021} with expert-based corrections and additions for keys that cannot be directly retrieved from \citet{fachagentur_wind_an_land_uberblick_2021}. The definition of all exclusions can be found in the supplementary data. 
\par
In scenarios~2, 2a, 2b, and 3 nation-wide exclusions are applied. Scenario S2 Expansive refers to typical buffers at the lower end of the federal state's legislation. Additionally, the construction of wind turbines in forests is allowed. Scenarios~2a and 2b build upon scenario~2, but restrict construction in protected landscapes and forests, respectively. Scenario~3 refers to exclusions on the higher end of the setback distances in the legislation. \autoref{tab:wind_exclusions} displays the most relevant of these.
\par
Inner areas are an influential exclusion, which refer to coherent built-up areas in accordance with §34~BauGB \cite{noauthor_baugesetzbuch_2017} and are modeled as the \textit{AX\_Ortslagen} features of Basis-DLM \cite{geobasisdaten__geobasis-de__bkg_2021_digitales_2021}. The residential buildings and other residential areas are therefore treated as areas without a development plan (outer areas). In the case of some the federal state laws, residential areas in outer areas with a statute (\textit{Außenbereichssatzung}) are protected identical manner to inner areas. For scenario 1, we assume that all residential buildings are protected by such a statute. 
\par
Following the land eligibility analysis, areas smaller than \SI{0.01}{\kilo\metre\squared} are excluded.
\begin{table*}
\captionsetup{justification=raggedright}
\caption{Selected exclusion criteria for the scenarios of the wind onshore potential analysis.}\label{tab:wind_exclusions}
\begin{tabular*}{\textwidth}{c|c||c|c|c|c|c}
\toprule
  Criterion & Data Source & S1\textsuperscript{1} & S2\textsuperscript{2} & S2a\textsuperscript{2a} & S2b\textsuperscript{2b} & S3\textsuperscript{3}  \\ 
\midrule
 Inner areas & Basis-DLM \cite{geobasisdaten__geobasis-de__bkg_2021_digitales_2021} & individual\textsuperscript{*} & \SI{1000}{\metre} & \SI{1000}{\metre} & \SI{1000}{\metre} & \SI{1000}{\metre} \\
 Residential buildings, outer areas & Hausumringe \cite{geobasisdaten__geobasis-de__bkg_2021_amtliche_2021}  & individual\textsuperscript{*} & \SI{3}{H} & \SI{3}{H} & \SI{3}{H} & \SI{1000}{\metre} \\
 Forests & Basis-DLM \cite{geobasisdaten__geobasis-de__bkg_2021_digitales_2021} & individual\textsuperscript{*} & not excluded & not excluded &\SI{0}{\metre}  & \SI{0}{\metre} \\
 Protected Landscapes & WDPA \citep{unep-wcmc_iucn_world_2016}  & individual\textsuperscript{*} & not excluded &  \SI{0}{\metre} & not excluded &  \SI{0}{\metre}  \\
\bottomrule
\end{tabular*}
\begin{tablenotes}
        \small{
        \item[1]\textsuperscript{1}S1 Legislation;\textsuperscript{2}S2 Expansive;\textsuperscript{2a}S2a No Protected Landscapes;\textsuperscript{2b}S2b No Forests;\textsuperscript{3}S3 Restrictive.}
        \item[2]\textsuperscript{*}federal state-specific exclusions
\end{tablenotes}
\end{table*}

To estimate the capacity potential, we classify turbines with \textit{GLAES}~\cite{ryberg_future_2018} in the eligible areas.
To this end, a spacing distance of \SI{8}{D}$\times$\SI{4}{D} is used. A typical low wind (wind class IEC IIIB) \SI{4.7}{\mega\watt} turbine with a \SI{155}{\metre} rotor diameter and \SI{120}{\metre} hub height is used as the reference turbine model.

\subsubsection{Results \& Discussion}\label{sssec:onwind-res}

\autoref{tab:wind_results} shows the main results of the scenarios on the national level. In scenario S2 Expansive, $7.25\%$ of Germany's area is usable for wind turbines, which leads to a capacity potential of \SI{404}{\giga\watt}. Scenario S3 Restrictive leads to only \SI{90}{\giga\watt} on $1.1\%$ of Germany's land. The capacity densities in the scenarios S1 Legislation, S2 Expansive, and S2a No Protected Landscapes are smaller than in scenarios S2b No Forests S2b and S3 Restrictive 
due to larger contiguous areas in the scenarios not excluding forests, i.e., S1, S2, and S2a. In the larger areas, the distance between the turbines is more apparent as smaller areas naturally promote spacing between them. 
\begin{table}
\caption{Results for the scenarios of the onshore wind potential analysis at the national level.}\label{tab:wind_results}
\begin{tabular*}{0.49\textwidth}{c||c|c|c|c|c}
\toprule
  {} & S1\textsuperscript{1} & S2\textsuperscript{2} & S2a\textsuperscript{2a} & S2b\textsuperscript{2b} & S3\textsuperscript{3}   \\ 
\midrule 
 Area [\si{\kilo\metre\squared}] & 24663 & 25938 & 17613 & 10056 & 3923 \\
 Area Share [\si{\%}] & 6.89 & 7.25 & 4.92 & 2.81 & 1.10 \\
 Capacity [\si{\giga\watt}] & 385 & 403 & 287 & 241 & 90 \\
 Density [\si[per-mode=fraction]{\mega\watt\per\kilo\metre\squared}] & 15.6 & 15.5 & 16.3 & 23.9 & 22.8 \\
\bottomrule
\end{tabular*}
\begin{tablenotes}
        \small{
        \item[1]\textsuperscript{1}S1 Legislation; \textsuperscript{2}S2 Expansive; \textsuperscript{2a}S2a No Protected Landscapes; \textsuperscript{2b}S2b No Forests; \textsuperscript{3}S3 Restrictive}
\end{tablenotes}
\end{table}
\par

Two potential analyses on the federal state level \cite{landesamt_fur_natur_umwelt_und_verbraucherschutz_nordrhein-westfalen_lanuv_potenzialstudie_2021,landesanstalt_fur_umwelt_baden-wurttemberg_potenzialanalyse_2019} are used to validate the workflow. The area potential adds up to $103.3\%$ and $96\%$ in relation to the study in Baden-Württemberg \cite{landesanstalt_fur_umwelt_baden-wurttemberg_potenzialanalyse_2019} and Northrhine-Westphalia \cite{landesamt_fur_natur_umwelt_und_verbraucherschutz_nordrhein-westfalen_lanuv_potenzialstudie_2021}, respectively. When compared with the \citet{landesanstalt_fur_umwelt_baden-wurttemberg_potenzialanalyse_2019}, an Intersection over Union of $79.3\%$ is achieved. Differences can be explained by the versions of Basis-DLM~\cite{geobasisdaten__geobasis-de__bkg_2021_digitales_2021} and unclear definitions of single exclusions.
\begin{figure}[h]
	\centering
\includegraphics[width=0.50\textwidth]{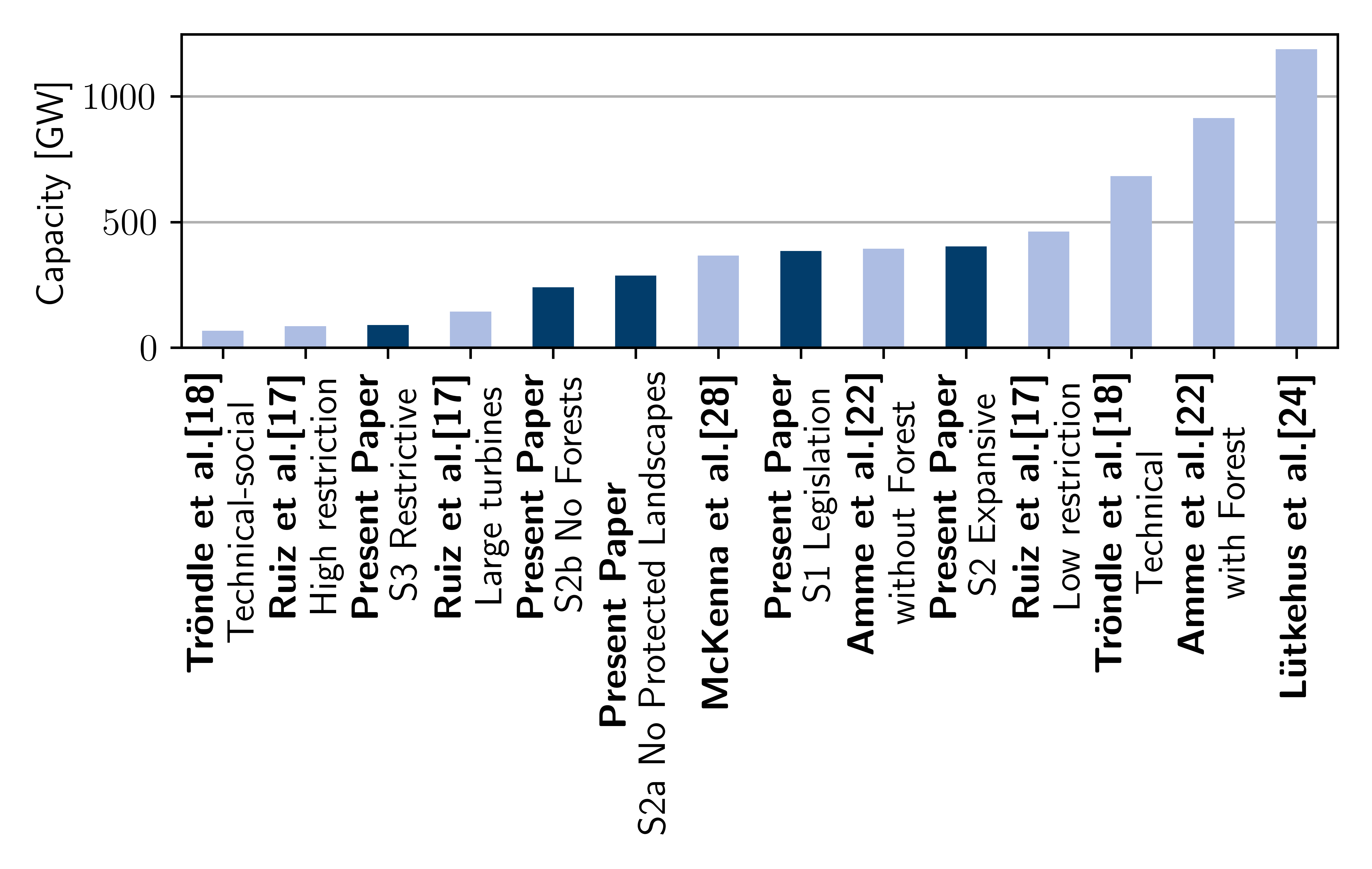}
	  \caption{Comparison of onshore wind potentials for studies providing capacity potential at the national level.}\label{fig:comparison_onshore}
\end{figure}
\par
\autoref{fig:comparison_onshore} shows the range of \SI{68}{\giga\watt} to \SI{1188}{\giga\watt} of onshore wind potential values in literature.
The results are influenced by the use of different datasets (see \autoref{sssec:data}), exclusion definitions, area correction factors, and capacity estimation methods.
\citet{trondle_home-made_2019,ruiz_enspreso_2019} and \citet{mckenna_cost-potential_2014} utilize, amongst others, CLC for land use classification, which we show to be problematic for potential analyses (\autoref{sssec:data}).
Furthermore, \citet{trondle_home-made_2019}, in technical-social-scenarios and \citet{mckenna_cost-potential_2014} correct the estimated areas by a correction and suitability factor for land use categories, respectively.
These approaches are not comparable to our results, which clearly indicate locations as being either eligible or ineligible.
\citet{amme_photovoltaik-_2021} outlined a scenario excluding residential areas with a \SI{1000}{\metre} setback distance as well as protected landscapes but allowing forests, which is comparable to our scenario S2a No Protected Landscapes.
By comparing the results, \citet{amme_photovoltaik-_2021} exceed the presented eligible area and capacity by a factor of 2.47 and 3.18, respectively.
\citet{lutkehus_potenzial_2013} use comparably low setback distances, e.g. \SI{600}{\metre} to inner areas, and neglect residential buildings, leading to the highest capacity potential of the considered studies. This emphasizes the impact of the chosen exclusion criteria and setback distances. Furthermore, they use the land use dataset DLM250~\citep{geobasisdaten__geobasis-de__bkg_2012_digitales_2012} with a positional accuracy of $\pm$\SI{100}{\metre}.
\par
There is no consensus on exclusions and setback distances in wind potential analyses. 
One exclusion used by all of the cited studies is residential areas. 
Nevertheless, individual residential buildings are often neglected due to data availability issues. 
In our scenarios, neglecting residential buildings leads to significantly higher capacity potentials ($27\%$, $9\%$, $10\%$, $12\%$, and $34\%$). 
For regional or federal state analyses, this error is even more apparent: For example, in Schleswig Holstein, the potential gets overestimated by $56\%$, $25\%$, $25\%$, $27\%$, and $100\%$. 
To the knowledge of the authors, no national potential analysis with published results has considered buildings. 
In future work, a focus on the sensitivity regarding exclusion definitions could help make such effects more transparent. 
Additionally, the sensitivities of certain exclusion criteria, e.g., forests or buildings, could help legislators quantify the impact of their decisions. 
Our restrictive scenario results in a potential of \SI{90}{\giga\watt}, which is not sufficient to reach Germany's climate goals according to different studies~\citep{stolten_strategien_2021,luderer_deutschland_2021,goke_100_2021,brandes_wege_2021} and underlies the importance of liberal legislation towards wind expansion.  
\par
Furthermore, the chosen turbine has a large impact on the capacity density, as well as on the land eligibility analyses via height- or diameter-dependent exclusions. We selected a typical wind turbine for low wind speeds (cf. \autoref{sssec:onwind-met}). When using a turbine for medium wind speeds (IEC II, \SI{5}{\mega\watt} capacity,  \SI{145}{\meter} diameter) with the same hub height, the capacity potential increases by $16.2\%$, to $19.5\%$, in our scenarios. 
In future, site-specific turbine selection could help improve the results in this regard. Nevertheless, with respect to the chosen low wind speed turbine, the results are still robust. Only $2\%$ and $29\%$ of the determined locations, respectively, reach an average wind speed  IEC II wind class in \SI{100}{\metre} and \SI{150}{\metre} height \cite{technical_university_of_denmark_global_2021}.

\subsection{Offshore wind potential}\label{ssec:offwind}
\subsubsection{Literature}\label{sssec:offwind-lit}
The following section provides an overview of offshore wind potential analyses to be found in the literature.
The regional coverage of the analyzed studies varies between global \citep{sensfus_langfristszenarien_2021,bosch_temporally_2018}, European \citep{ebner_regionalized_2019,trondle_home-made_2019,ruiz_enspreso_2019,zappa_analysing_2018} and national \citep{luderer_deutschland_2021}. 
As a first step in the potential analyses, eligible areas are determined. 
Most studies \citep{bosch_temporally_2018,ebner_regionalized_2019,trondle_home-made_2019, caglayan_techno-economic_2019,zappa_analysing_2018,ruiz_enspreso_2019} use greenfield approaches, initially considering the sea to be eligible and applying exclusion criteria.
\citet{sensfus_langfristszenarien_2021-1}, however, follow a mixed approach of greenfield analysis and pre-selected areas in their global analysis. 
For Germany, \SI{80}{\%} of the declared offshore wind areas in 4C Offshore~\citep{4c_offshore_global_2021} are used as pre-selected areas.
\citet{luderer_deutschland_2021} consider the wind farm areas of the draft of the Area Development Plan 2020 \citep{bundesamt_fur_seeschifffahrt_und_hydrographie_entwurf_2020} and areas of existing plants \citep{directorate-general_for_maritime_affairs_and_fisheries_european_2021} as pre-selected ones. 

The studies that follow the greenfield approach often use similar exclusion criteria, but differ in terms of dataset, methods, and buffer distances. 
Commonly used exclusion criteria include, amongst others, distance to shore, sea depth, protected areas, shipping routes, and infrastructures.
The minimal distance to shore varies between \SI{10}{\kilo\metre}~\citep{bosch_temporally_2018} and \SI{22.2}{\kilo\metre}~\citep{ruiz_enspreso_2019}.
The exclusion of areas with large sea depths differs between depths lower than
\SI{50}{\metre}~\citep{trondle_home-made_2019,ruiz_enspreso_2019,zappa_analysing_2018}, \SI{100}{\metre}~\citep{ruiz_enspreso_2019}, \SI{1000}{\metre}~\citep{bosch_temporally_2018,ebner_regionalized_2019,caglayan_techno-economic_2019}, or unlimited \citep{zappa_analysing_2018}.
Shipping routes are excluded by \citet{ruiz_enspreso_2019} and \citet{caglayan_techno-economic_2019} based on the dataset of \citet{halpern_cumulative_2015}, which addressed the likelihood of ships.
Furthermore, \citet{caglayan_techno-economic_2019} manually exclude known routes with a buffer of \SI{4}{\kilo\metre}.
\citet{ebner_regionalized_2019} exclude shipping routes, but without noting a methodology or data source.
Additionally, certain studies \citep{caglayan_techno-economic_2019,ruiz_enspreso_2019,bosch_temporally_2018} exclude infrastructure, e.g., cables and pipelines with buffers from \SI{500}{\metre} \citep{caglayan_techno-economic_2019} to \SI{7.4}{\kilo\metre} \citep{ruiz_enspreso_2019}.
After applying the exclusion criteria, \citet{zappa_analysing_2018} and \citet{trondle_home-made_2019}, in their technical-social, scenario reduce the resulting eligible areas to \SI{20}{\%} and \SI{10}{\%}, respectively.
\citet{zappa_analysing_2018} explain the reduction by a set factor due to missing exclusion categories.

The potential offshore wind capacity can be estimated for the eligible areas.
Most presented studies employ an aggregated method with capacity density factors ranging from \SI{3.14}{\mega\watt\per\kilo\metre\squared}~\citep{bosch_temporally_2018} to \SI{15}{\mega\watt\per\kilo\metre\squared}~\citep{trondle_home-made_2019}.
\citet{ebner_regionalized_2019} reduce \SI{14}{\mega\watt\per\kilo\metre\squared} to \SI{5}{\mega\watt\per\kilo\metre\squared} due to a results deviation of a factor of 3 to the German Bundesfachplan 2017~\citep{bundesamt_fur_seeschifffahrt_und_hydrographie_bundesfachplan_2017,bundesamt_fur_seeschifffahrt_und_hydrographie_bundesfachplan_2017-1}.
Other studies~\citep{sensfus_langfristszenarien_2021-1,luderer_deutschland_2021} use a variable capacity density per eligible area.
Only \citet{caglayan_techno-economic_2019} employ a placing algorithm for wind turbines with a turbine spacing of \SI{10}{D}$\times$\SI{4}{D}. 

\subsubsection{Methodology}\label{sssec:offwind-met}
The offshore wind potential analysis was performed for  federal states in the coastal sea and Exclusive Economic Zones (EEZ) in the Northern and Baltic sea areas.
Four different scenarios were considered: 

\begin{itemize}
    \item S1 Expansive: Greenfield analyses with offshore wind expansion favoring exclusions \\ S1a  Military: S1, including usage of military areas
    \item S2 Legislation: Current priority and reservation areas for offshore wind in legislation
    \item S3 Restrictive Legislation: Current priority areas for offshore wind in legislation
\end{itemize}

The Low Exclusion Scenarios S1 and S1a are greenfield approaches that use exclusion criteria with comparably low exclusion definitions. 
The following describes the main exclusions; full information about these can be found in the supplementary data.
Areas within \SI{15}{\kilo\metre} of shore and \SI{500}{\metre} to neighboring seas are considered ineligible.
To address shipping land use, the declared priority shipping areas of the current legislation of the EEZ~\citep{bundesamt_fur_seeschifffahrt_und_hydrographie_raumordnungsplan_2021} and federal states~\citep{niedersachsisches_ministerium_fur_ernahrung_landwirtschaft_und_vebraucherschutz_neubekanntmachung_2017,ministerium_fur_energie_infrastruktur_und_digitalisierung_landesraumentwicklungsprogramm_2016,landesplanung_schleswig-holsteinmilig_fortschreibung_2021} are excluded with a buffer of \SI{500}{\metre}.
Furthermore, infrastructure facilities, e.g., cables and platforms, are excluded by CONTIS~\citep{bundesamt_fur_seeschifffahrt_und_hydrographie_contis_2020} with \SI{500}{\metre}.
The scenarios S1 Expansive and S1a Military differ in the consideration of designated military areas, which cover a significant share of the German coastlines. 
As military areas overlap with designated wind areas, a complete exclusion is unrealistic. 
However, the exclusion of individual sub-areas is not possible due to a lack of indications of the suitability of mixed use.
Only Scenario S1 Expansive further excludes declared military areas according to CONTIS~\citep{bundesamt_fur_seeschifffahrt_und_hydrographie_contis_2019} and ROP~\citep{bundesamt_fur_seeschifffahrt_und_hydrographie_raumordnungsplan_2021}.

The scenario S2 Legislation  considers all designated areas for offshore wind as pre-selected, whereas scenario S3 Restrictive Legislation  only considers the priority and conditional priority areas ones. 
The definition of these areas is based on the current legislation of the EEZ~\citep{bundesamt_fur_seeschifffahrt_und_hydrographie_raumordnungsplan_2021}, Lower Saxony~\citep{niedersachsisches_ministerium_fur_ernahrung_landwirtschaft_und_vebraucherschutz_neubekanntmachung_2017}, and Mecklenburg-Western Pomerania~\citep{ministerium_fur_energie_infrastruktur_und_digitalisierung_landesraumentwicklungsprogramm_2016}. 
Schleswig Holstein and Hamburg do not designate areas for offshore wind~\citep{bundesamt_fur_seeschifffahrt_und_hydrographie_flachenentwicklungsplan_2020,landesbetrieb_geoinformation_und_vermessung_lgv_hamburg_flachennutzungsplan_2021}.
The Legislation scenarios S2 and S3 further apply the infrastructure exclusions described in Low Exclusion Scenarios S1 and S1a. 

After employing the exclusions of the scenarios, areas smaller than \SI{0.1}{\kilo\metre\squared} are excluded.

The potential capacity of the scenarios is estimated with a similar approach as that used for onshore wind (\autoref{sssec:onwind-met}) by using the turbine placement of \textit{GLAES}~\cite{ryberg_future_2018} and a turbine spacing of \SI{8}{D}$\times$\SI{4}{D}, respectively.
A reference turbine with a capacity of \SI{8}{\mega\watt} and a rotor diameter of \SI{167}{\metre} is used.

\subsubsection{Results \& Discussion}\label{sssec:offwind-res}

The results of the offshore potential analysis on the national level are presented in \autoref{tab:offshore_wind_results} and show a range of between \SI{34.1}{\giga\watt} (S3 Restrictive Legislation) and \SI{99.6}{\giga\watt} (S1a Military). 

\begin{table}
\caption{Results for the scenarios of the offshore wind potential analysis on a national level}\label{tab:offshore_wind_results}
\begin{tabular*}{0.42\textwidth}{c||c|c|c|c}
\toprule
  {} & S1\textsuperscript{1} & S1a\textsuperscript{1a} & S2\textsuperscript{2} & S3\textsuperscript{3}
  \\ 
\midrule 
 Area [\si{\kilo\metre\squared}] & 7353 & 9275 & 5174 &  3182 \\
 Area Share [\si{\%}] & 13.07 & 16.48 & 9.19 & 5.65  \\
 Capacity [\si{\giga\watt}] & 79.1 & 99.6 & 55.8 & 34.1  \\
 Density [\si[per-mode=fraction]{\mega\watt\per\kilo\metre\squared}] & 10.75 & 10.74 & 10.79 & 10.71   \\
\bottomrule
\end{tabular*}
\begin{tablenotes}
        \small{
        \item[1]\textsuperscript{1}S1 Expansive; \textsuperscript{1a}S1a Military; \textsuperscript{2}S2 Legislation; and \textsuperscript{3}S3 Restrictive Legislation}
\end{tablenotes}
\end{table}

It can be noted that there exists the possibility to increase the potential from current legislation, i.e. S2 and S3, by limiting other declared land uses such as shipping  or military areas.
However, the possible increase in capacity potential highly differs between the Northern and  Baltic seas. 
From scenario S2 Legislation to S1 Expansive the potential is increased by \SI{21.6}{\giga\watt} in the North Sea and \SI{1.6}{\giga\watt} in Baltic.
For the comparison of S2 Legislation to S1a Military, the capacity increase is higher, with \SI{39}{\giga\watt} and \SI{5}{\giga\watt}, respectively.
Thus, by redefining and limiting other declared land uses in legislation as military or shipping areas, further areas could be designated for offshore wind, especially in the North Sea.
Future work could critically revise these areas or consider combined usage with offshore wind projects.

Furthermore, the reference turbine impacts the capacity potential.
In the future, site-specific turbine designs and a sensitivity analysis for turbine selection could be applied. 

\begin{figure}[h]
	\centering
\includegraphics[width=0.50\textwidth]{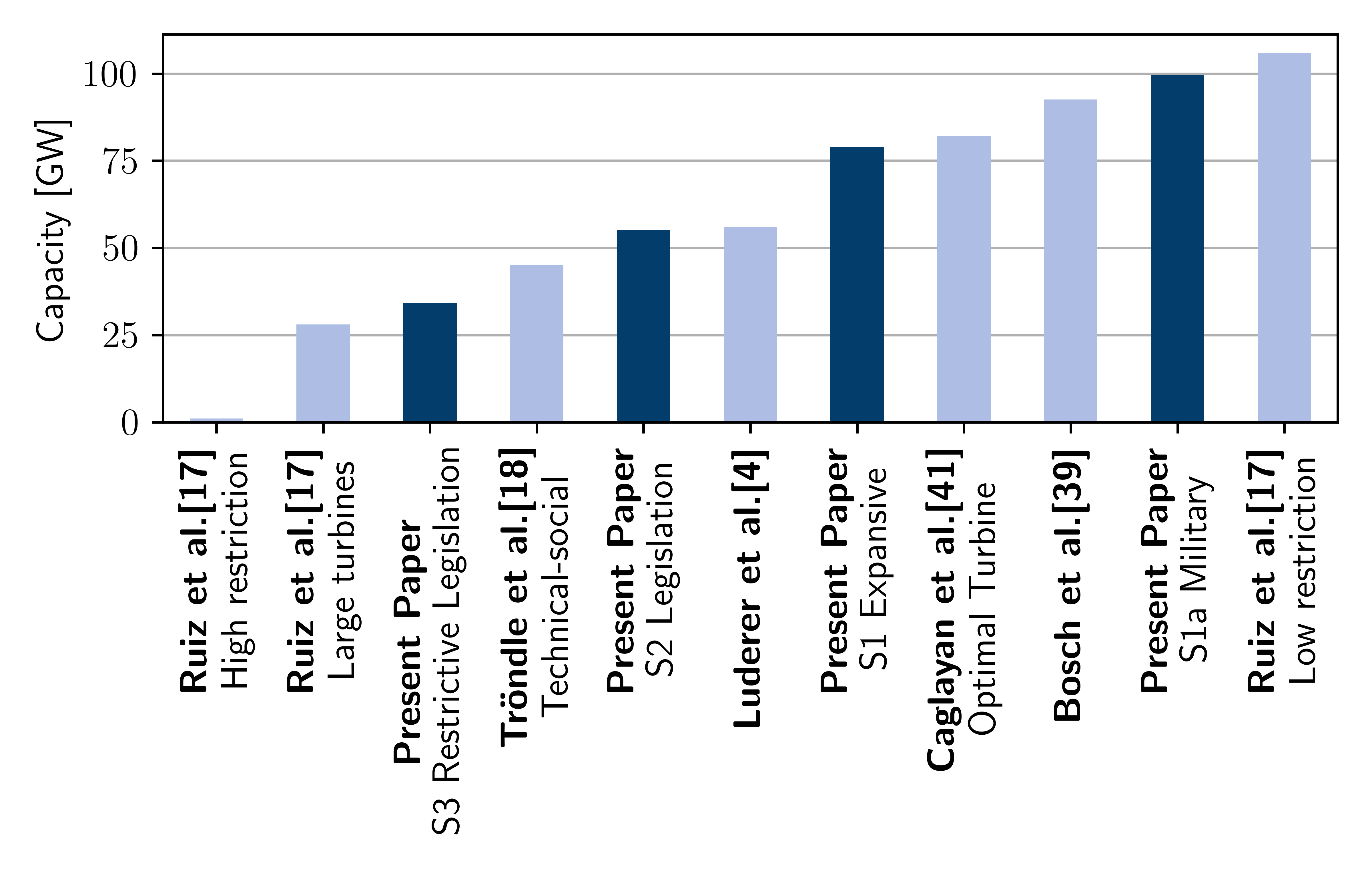}
	  \caption{Comparison of offshore wind potentials for studies providing capacity potential on the national level}\label{fig:comparison_offshore}
\end{figure}

The results of the presented scenarios are compared to the literature in \autoref{fig:comparison_offshore}. 
The scenarios S1 Expansive and S1a Military follow a greenfield approach, comparable to all considered studies except for that of \citet{luderer_deutschland_2021}.
\citet{luderer_deutschland_2021}, who use the same designated wind farm area of BSH leading to \SI{56}{\giga\watt}, comparable to \SI{55.8}{\giga\watt} in the present scenario, S2 Legislation. 

For the greenfield analyses, none of the presented studies use official German data, e.g., for shipping or military areas.
Two scenarios of \citet{ruiz_enspreso_2019} and the technical-social scenario of \citet{trondle_home-made_2019} even undercut the capacity potential results of our scenarios S2 Legislation and S3 Restrictive Legislation.
\citet{ruiz_enspreso_2019} in their low exclusion scenario and \citet{bosch_temporally_2018} have comparable capacity potentials to S1a.
However, as the capacity densities of our study are around two and three times higher in comparison to their studies, respectively, this indicates significantly fewer eligible areas in our scenarios.
The comparisons show the impact of data sources, exclusions, and capacity densities on the land eligibility analyses and capacity estimations.

\subsection{Open-field photovoltaic potential}\label{ssec:ofpv}

\subsubsection{Literature}\label{sssec:ofpv-lit}
Although wind potential analyses mostly use greenfield approaches in land eligibility analysis~(see \autoref{ssec:onwind}), open-field PV potential analyses often first consider pre-selected areas as being eligible and then apply further exclusion criteria \citep{lux_langfristszenarien_2021,ruiz_enspreso_2019,trondle_home-made_2019,ebner_regionalized_2019,landesanstalt_fur_umwelt_baden-wurttemberg_potenzialanalyse_2018,seidenstucker_solarkataster_2020,amme_photovoltaik-_2021,luderer_deutschland_2021}. 
Only a few potential analyses choose a greenfield exclusion approach \citep{peters_raumlich_2015,ryberg_generation_2019}.

However, as the pre-selected areas highly differ accross the studies, several examples are presented: \citet{ruiz_enspreso_2019} and \citet{trondle_home-made_2019} estimate the potential for Europe and 
consider land use categories as pre-selected areas, which leads to a high area potential for open-field PV.
\citet{ruiz_enspreso_2019} define amongst others, bare land, some cropland, and some eligible classes of vegetation and use a combination of Global Land Cover~\citep{european_environment_agency_global_2016} and Corine Land Cover~\citep{copernicus_programme_corine_2018} to identify these. 
This leads to an initial eligible area of $44\%$ of Germany, from which a share of $3\%$ is considered eligible in another scenario. 
\citet{trondle_home-made_2019} consider $10\%$ of bare and unused land defined by GlobeCover2009~\citep{esa_globcover_2009} as pre-selected areas in a technical-social potential analysis. 
\citet{lux_langfristszenarien_2021} define certain shares of categories of CLC~\citep{copernicus_programme_corine_2018} eligible, e.g., $16\%$ of bare land and $2\%$ of bushland.
Meanwhile, \citet{ebner_regionalized_2019} consider agriculture and grazing land identified by CLC \citep{copernicus_programme_corine_2018} in less favored regions \citep{european_environment_agency_less_2012} eligible for technical potential.
However, this is further reduced to $7\%$ of the pre-selected area, which is $50\%$ of the area share currently used for energy crops.
Several studies with the scope of Germany follow the legislation of the Renewable Energy Act (EEG)~\citep{noauthor_gesetz_2021}. 
Amongst other areas, the EEG considers less-favored areas and side strips of motorways and railways to be eligible. 
The potential areas along side strips are either defined by a \SI{110}{\metre} width in accordance with EEG~2017~\citep{noauthor_gesetz_2017} or a \SI{200}{\metre} width in accordance with EEG 2021~\citep{noauthor_gesetz_2021} minus a \SI{15}{\metre} animal migration buffer.
\citet{luderer_deutschland_2021} also use side strips, but consider \SI{185}{\metre} around the line objects of OSM, therefore neglecting the width of the road and the animal migration corridor.
Furthermore, they consider agricultural land with poor soil quality, which is defined with the dataset Soil Quality Rating~(SQR)~\citep{bundesanstalt_fur_geowissenschaften_und_rohstoffe_bgr_ackerbauliches_2013} and a threshold of 40.
\citet{amme_photovoltaik-_2021} similarly consider the less-favored areas, but also consider side strips with a width of \SI{500}{\metre}, which corresponds to 2.7 times the current legislation. Agricultural land with higher soil quality rating than 40 are excluded within the side strips. 
Other EEG-favored areas are also neglected due to data availability and small resulting areas.
Several studies perform potential analyses on a federal state level in Germany, e.g., \citet{seidenstucker_solarkataster_2020} for North Rhine-Westphalia and \citet{landesanstalt_fur_umwelt_baden-wurttemberg_potenzialanalyse_2018} for Baden-Württemberg. 
Both studies use input datasets with local coverage and high resolutions, e.g., Basis-DLM \cite{geobasisdaten__geobasis-de__bkg_2021_digitales_2021}, and estimate the potential for the side strips of \SI{110}{\metre} of rails and road ways in accordance with the legislation of EEG 2017 \citep{noauthor_gesetz_2017}. 
\par
After defining the pre-selected areas, further exclusion criteria can be applied, resulting in eligible areas of various shapes and sizes.
As it is not economically-feasible to install open-field PV sites in small eligible areas, several studies exclude areas smaller than a certain threshold. 
The threshold can vary in the range from \SI{500}{\metre\squared} \citep{seidenstucker_solarkataster_2020} to \SI{100,000}{\metre\squared} \citep{amme_photovoltaik-_2021}.

As the second step of the potential analysis, the installable open-field PV capacity on eligible areas is estimated.
Open-field PV potential analysis uses an empirical factor of the capacity density in \SI{}{\mega\watt\per\km\squared}. 
However, the literature reports a high range from \SI{40}{\mega\watt\per\km\squared} in \citep{ebner_regionalized_2019} to \SI{300}{\mega\watt\per\km\squared} \citep{ruiz_enspreso_2019}.
\citet{ebner_regionalized_2019} explain the factor of \SI{40}{\mega\watt/\km\squared} based on aerial photo evaluation of existing open-field PV systems. 
\citet{trondle_home-made_2019} explain their value of \SI{80}{\mega \watt \per \km\squared} with a module efficiency of $16\%$ and $50\%$ area reduction by row placements to prevent shadowing.  
\citet{seidenstucker_solarkataster_2020} assume a module efficiency of $17\%$ and also an area reduction of $50\%$ for slopes less than 20°, leading to \SI{85}{\mega\watt/\km\squared}.
\citet{wirth_recent_2021} reports a value of \SI{100}{\mega\watt\per\km\squared} based on row spacing and a module efficiency of $20\%$.

\subsubsection{Methodology}\label{sssec:ofpv-met}
In the first step of the potential analysis for open-field PV pre-selected areas are considered eligible followed by further exclusions. 
To this end three scenarios are regarded: 

\begin{itemize}
    \item S1 Side Strips: Side strips of motorways and railways in accordance with the subsidy areas of the EEG 2021 \citep{noauthor_gesetz_2021}
    \item S2 Poor Soil: Arable land with barren soil based on the Soil Quality Rating~(SQR)~of~\citep{bundesanstalt_fur_geowissenschaften_und_rohstoffe_bgr_ackerbauliches_2013} 
    \item S3 Combination: Side strips, for which arable land is restricted to SQR $<40$, and S2 
\end{itemize}

Scenario S1 Side Strips estimates the potential of the side strips of motorways and railways. 
Their routes and widths are identified by Basis-DLM~\cite{geobasisdaten__geobasis-de__bkg_2021_digitales_2021}.
Then, side strips of \SI{200}{\metre} are used, from which \SI{15}{\metre} on the inside is subtracted for animal migration. 
The remaining area represents the pre-selected area.

\begin{figure}[h]
	\centering
\includegraphics[width=0.49\textwidth]{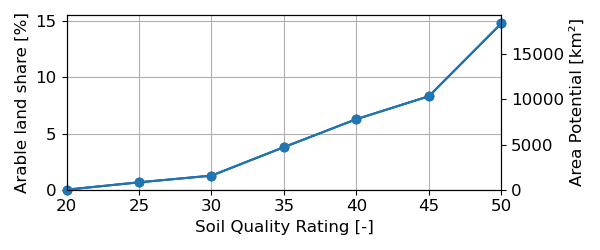}
	  \caption{Sensitivity analysis of the Soil Quality Rating threshold for arable land in the land eligibility analysis for open-field PV.}\label{fig:SA_SQR}
\end{figure}

Scenario S2 Poor Soil pre-selects arable land with a bad soil quality based on the SQR dataset~\citep{bundesanstalt_fur_geowissenschaften_und_rohstoffe_bgr_ackerbauliches_2013}.  
The SQR threshold was chosen to be 30 based on a pre-analysis.
\autoref{fig:SA_SQR} shows the area potential and share of arable land for SQR limits of between 20 and 50.
As is shown in \autoref{sssec:ofpv-lit}, other potential studies use a threshold of 40, which results in $6.28\%$ of arable land.
We assume this to be too ambitious due to land use conflicts with cultivation of arable land and therefore choose the threshold of 30.
All areas up to a selected SQR value are intersected with the arable land in Germany from Basis-DLM \cite{geobasisdaten__geobasis-de__bkg_2021_digitales_2021} due to the resolution of \SI{100}{\metre}$\times$\SI{100}{\metre} in the BGR rating \cite{bundesanstalt_fur_geowissenschaften_und_rohstoffe_bgr_ackerbauliches_2013}. 

\par
Scenario S3 Combination is a combination of a pre-selection using side strips and arable land with barren soil. 
As side strips are considered less valuable than other areas, arable land with SQR values of less than 40 are excluded. 
On all other arable land, only SQR values of less than 30 are used as a pre-selection. 
It should be noted that other land uses than agricultural inside the side strips are not restricted by the SQR rating.

\par

After identifying pre-selected eligible areas for each scenario, unsuitable areas within these are deducted by further exclusion criteria.
\autoref{tab:ofpv_exclusions} displays the most relevant exclusions. 
The definition of all exclusions can be found in the supplementary material. 
After the land eligibility analysis, eligible areas smaller than a threshold of \SI{5000}{\metre\squared} were excluded for economic reasons, which lies within the range of the presented literature (see \autoref{sssec:ofpv-lit}).
However, this area threshold leads to a deviation of eligible areas at the different regional levels. 
Areas may become ineligible if split by a border to a size lower than the threshold, whereas on a higher regional level still eligible.
For the three scenarios, this leads to a capacity deviation from $1.3\%$ to $1.5\%$ between the federal state and the municipality level.

\par
\begin{table*}
\captionsetup{justification=raggedright}
\caption{Selected exclusion criteria for the scenarios of the open-field PV potential analysis.}\label{tab:ofpv_exclusions}
\begin{tabular*}{\textwidth}{c|c||c|c|c}
\toprule
  Criterion & Data Source & S1 Side Strips & S2 Poor Soil & S3 Combination   \\ 
\midrule
 Forests & Basis-DLM \cite{geobasisdaten__geobasis-de__bkg_2021_digitales_2021} & \SI{10}{\metre} & \SI{10}{\metre} & \SI{10}{\metre} \\
 All Buildings & Hausumringe \cite{geobasisdaten__geobasis-de__bkg_2021_amtliche_2021}  & \SI{10}{\metre} & \SI{10}{\metre} & \SI{10}{\metre} \\
  Arable land & Basis-DLM \cite{geobasisdaten__geobasis-de__bkg_2021_digitales_2021}, SQR \cite{bundesanstalt_fur_geowissenschaften_und_rohstoffe_bgr_ackerbauliches_2013}  & not excluded & SQR$\geq30$ & SQR$\geq30$, Sidestripes: SQR$\geq40$  \\
 Motorways, Railways & Basis-DLM \cite{geobasisdaten__geobasis-de__bkg_2021_digitales_2021}  & \SI{15}{\metre} & \SI{200}{\metre} & \SI{15}{\metre}  \\
\bottomrule
\end{tabular*}
\end{table*}

As a next step in the potential analysis, the installable capacity on the eligible land was estimated.
Therefore, the eligible areas and their size $A_{OFPV}$ are extracted, on which the modules are placed in a southerly direction with an optimal tilt for the latitude, which is calculated using the \textit{RESKit} tool~\cite{ryberg_evaluating_2018} in accordance with \citet{ryberg_generation_2019}.
This area can be converted by taking a factor for row spacing $f_{RS}=0.5$, a factor for construction-related obstructions, e.g., access roads or surrounding areas, $f_{c}=0.72$ \cite{ong_land-use_2013}, and the efficiency of the modules, which is assumed to be $\eta_{PV}=0.22$ in line with current high-end modules \citep{suntech_ultra_2021,lg_lg_2021}, into consideration:

\begin{equation}
\begin{split}
    P_{OFPV} & = \eta_{OFPV}\cdot f_{RS}\cdot f_{c} \cdot A_{OFPV} \\
    & =  \SI{79.2}{\mega\watt\per\kilo\metre\squared} \cdot A_{OFPV}
\end{split}
\end{equation}

\subsubsection{Results \& Discussion}\label{sssec:ofpv-res}

The results at the national level for the potential area and the capacity of the scenarios are shown in \autoref{tab:ofpv_results}.
The capacity potential varies between \SI{123.6}{\giga\watt_p} in S2 Poor Soil and \SI{456.1}{\giga\watt_p} in S1 Side Strips.

\begin{table}
\caption{Results for the scenarios of the open-field PV potential analysis at the national level.}\label{tab:ofpv_results}
\begin{tabular*}{0.49\textwidth}{c||c|c|c}
\toprule
  {} & S1\textsuperscript{1} & S2\textsuperscript{2} & S3\textsuperscript{3}   \\ 
\midrule 
 Area [\si{\kilo\metre\squared}] & 5723 & 1560 & 4373  \\
 Area Share [\si{\%}] & 1.60 & 0.44 & 1.22 \\
 Capacity [\si{\giga\watt_p}] & 456.1 & 123.6 & 347.7  \\
\midrule 
 No. of Municipalities  & 11003 & 11003 & 11003 \\
 ...with pre-selected areas &  5667 & 1892 & 6446 \\
 ...with potential &  5253 & 1711 & 5939 \\
\bottomrule
\end{tabular*}
\begin{tablenotes}
        \small{
        \item[1]\textsuperscript{1}S1 Side Strips; \textsuperscript{2}S2 Poor Soil; and \textsuperscript{3}S3 Combination}.
\end{tablenotes}

\end{table}

To validate the presented workflow, two potential analyses for the federal state level
\citep{landesanstalt_fur_umwelt_baden-wurttemberg_potenzialanalyse_2018,seidenstucker_solarkataster_2020} were used due to their high spatial resolutions and usage of local data. 
Similarly to Side Strips S1, the studies consider pre-selected areas at the sides of roads and railways.
Using similar exclusion criteria and datasets as the corresponding studies, an Intersection over Union of $81.38\%$ for \citet{landesanstalt_fur_umwelt_baden-wurttemberg_potenzialanalyse_2018} and $87.88\%$ for \citet{seidenstucker_solarkataster_2020} were achieved. 
The area potential of the presented workflow in relation to that reported by  \cite{landesanstalt_fur_umwelt_baden-wurttemberg_potenzialanalyse_2018} and \cite{seidenstucker_solarkataster_2020} is $88\%$ and $103\%$, respectively. 
Remaining differences between the results of the studies and this paper arise due to different versions of land cover datasets and unclear definitions of single exclusions.
\par
\begin{figure}[h]
	\centering
\includegraphics[width=0.50\textwidth]{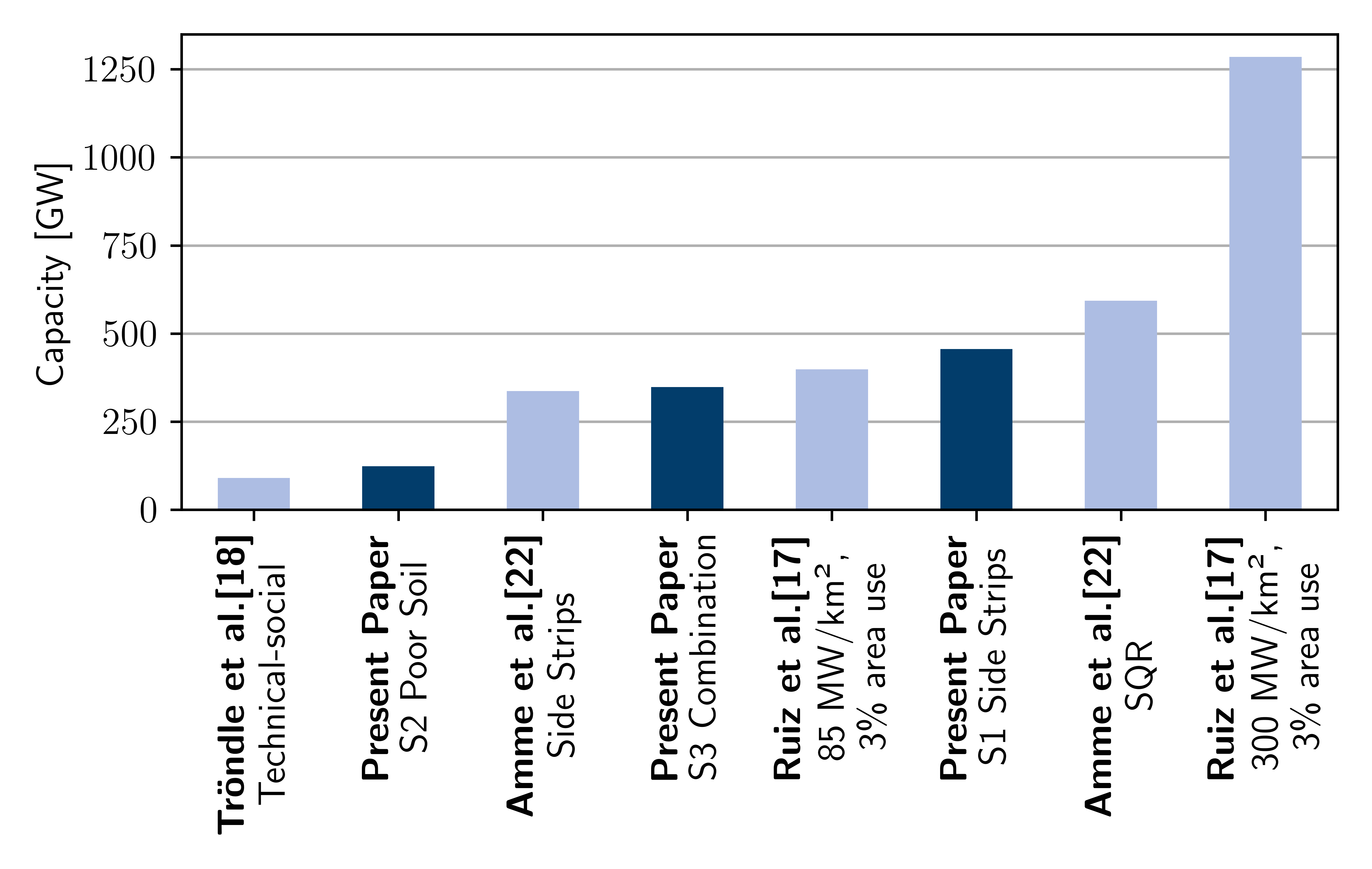}
	  \caption{Comparison of open-field PV potentials for studies providing capacity potential at the national level.}\label{fig:comparison_ofpv}
\end{figure}

As is shown in \autoref{fig:comparison_ofpv}, the literature provides a high range of \SI{90}{\giga\watt_p}~\citep{trondle_home-made_2019} to \SI{1285}{\giga\watt_p}~\citep{ruiz_enspreso_2019} as capacity potentials for open-field PV, because the pre-selected areas, used datasets, area reduction factors, and applied exclusion criteria highly differ between the studies (cf. \autoref{sssec:ofpv-lit}). 
Furthermore, the used capacity density factors vary greatly between the studies, which leads to further deviations.
\citet{ruiz_enspreso_2019} and \citet{trondle_home-made_2019} use set shares, $3\%$ and $10\%$, respectively, of large land use categories for the land eligibility analyses, and are therefore not comparable to our scenarios.

Methodology-wise our presented scenarios are only comparable to \citet{amme_photovoltaik-_2021}.
For the presented SQR scenario, the higher potential capacity of \citet{amme_photovoltaik-_2021} can be explained by its higher SQR threshold. 
For the side strips scenario, \citet{amme_photovoltaik-_2021} use side strips that are 2.7 times the width of current legislation. Nevertheless, by additionally excluding areas with a SQR higher than 40, the resulting capacity is lower than in our scenario S1 Side Strips.

\par

The presented workflow could be further extended by considering other subsidy-areas of the Renewable Energy Act~\citep{noauthor_gesetz_2021} like dumps and landfills as pre-selected areas.
Furthermore, categories beyond the scope of current legislation could be chosen as pre-selected areas. 
However, there exists no consensus in the literature about eligible land use categories and considering entire categories eligible tend to lead to high potentials, as in \citet{ruiz_enspreso_2019}. 
Such high area shares are critical due to the land use conflicts they can induce: Whereas wind turbines only occupy a share of their designated areas, open-field PV plants cover most of their appointed ones, which can lead to conflicts with, for example, agricultural land use. Furthermore, legislation has a significant impact on the distribution of potentials. Especially in the open-field PV scenarios, many municipalities have no potential (cf. \autoref{tab:ofpv_results}), due to the limitation to municipalities with subsidy-eligible areas. This can lead to unequally-distributed economic benefits among municipalities, as well as higher resistance in society because of visual impacts, for instance.
As the results are also highly sensitive towards the chosen capacity density factor, one avenue of further work could be to incorporate future system designs. Additionally, the effect of placing the modules flat or in an east-west orientation without any row-spacing could be analyzed.

\section{Rooftop photovoltaic potential}\label{sec:rtpv}
In the following chapter, the assessment of rooftop PV potential for Germany is discussed. First, an overview of the literature is given (\autoref{sssec:rtpv-lit}). Then, the methodology to extract potential from 3D building data is described (\autoref{sssec:rtpv-met}). The results are then presented and discussed (\autoref{sssec:rtpv-res}).

\subsection{Literature}\label{sssec:rtpv-lit}

Rooftop PV potential analyses can be classified by their geographical extent and the resolution of their results \citep{castellanos_rooftop_2017,mainzer_high-resolution_2014}.
\par
The resolution of different approaches can be split into high, medium, and low levels. 
High level approaches, which can assess the potential for individual buildings include, for example, image recognition techniques \cite{mainzer_assessment_2017,song_approach_2018,sampath_estimation_2019,singh_estimation_2013} or workflows using 3D building models \citep{grothues_landesweite_2018,strzalka_large_2012,walch_big_2020,fath_method_2015,jakubiec_method_2013}.
For instance, \citet{grothues_landesweite_2018} determine all roof geometries in North Rhine-Westphalia on the basis of laser scan data at \SI{0.5}{\metre}$\times$\SI{0.5}{\metre} resolutions. 
Low- and medium-level approaches base their analysis mostly on statistical data and are therefore limited to the assessment of potentials for the regional scope of the data, e.g., to regions, countries, or raster resolution \citep{ruiz_enspreso_2019,trondle_home-made_2019,peters_raumlich_2015,mainzer_high-resolution_2014,kurdgelashvili_estimating_2016,assouline_quantifying_2017,bodis_high-resolution_2019,ebner_regionalized_2019,schmid_regionalisierung_2021}. 
\citet{trondle_home-made_2019}, for instance, use population data from the European Settlement Map \cite{european_commission_joint_research_centre_european_2016} to estimate roof areas and calibrate them with data from sonnendach.ch~\citep{portmann_sonnendachch_2019}. 
Meanwhile, \citet{ruiz_enspreso_2019} use CLC~\citep{copernicus_programme_corine_2018} data to identify residential and industrial areas in which they estimate roof areas based on fixed factors.
\par
The geographical extent of these analyses varies from single buildings to continent-wide.
Statistical approaches are mainly used to perform large-scale analyses, such as at national \citep{peters_raumlich_2015,mainzer_high-resolution_2014,ebner_regionalized_2019,schmid_regionalisierung_2021} or continental scales \citep{ruiz_enspreso_2019,trondle_home-made_2019,bodis_high-resolution_2019}.
However, more recent studies were able to apply high-level methodologies to the national scale.
\citet{walch_big_2020} use LoD2 3D data to extract rooftop PV potential in Switzerland. 
To this end, they calculate the area and orientations of 9.6 million rooftops directly from the data. 
\citet{luderer_deutschland_2021} estimate the rooftop PV potential in Germany at the municipality level by correlating it with building footprints in a region. 
Similarly, \citet{wiehe_nothing_2021} used building footprints in Germany and building use to estimate generation potential. 
\citet{eggers_pv-ausbauerfordernisse_2020} use LoD1 building models (without information on rooftop geometry) to estimate German rooftop PV potential.

Approaches, that use data without information about superstructures, e.g., chimneys and windows, on roofs to estimate the available roof area, utilize reduction factors to adapt the potentially usable areas.
\citet{mainzer_high-resolution_2014} use a factor of $0.58$ to exclude obstacles on roofs and in areas with too much shadowing. In turn, \citet{fath_method_2015} differentiates between flat ($0.7$) and tilted ($0.75$) roofs based on the work of \citet{kaltschmitt_potentiale_1992}. 
The \citet{international_energy_agency_iea_potential_2002} determines a factor of 0.6 for constructions, shading, and historical elements. 
\citet{walch_big_2020} estimate the factor from LoD4 data in Geneva (38,000 roofs) with a machine learning approach for other LoD2 building models and estimate the factor to be between 0 and 0.8 depending on the roof size and tilt.
\citet{portmann_sonnendachch_2019} distinguish between flat ($0.7$) and tilted ($0.42-0.8$) roofs, for which a differentiation between roof size and category is carried out. 
\citet{eggers_pv-ausbauerfordernisse_2020} use a factor of $0.486$ to account for obstructions, shading, and inefficiencies when placing the modules. 
\citet{wiehe_nothing_2021} assume that $60\%$ of the roof area on residential and $80\$$ of the roof area on industrial and commercial structures are usable. 
\citet{grothues_landesweite_2018} estimate the impact of obstructions directly based on a highly resolved (\SI{0.5}{\metre}$\times$\SI{0.5}{\metre}) digital elevation model.
\par
In contrast to land eligibility analyses (see \autoref{sec:LEAPotentials}), in which a clear workflow has been established, the methodologies between different rooftop PV potential analyses vary greatly. Many high-resolution approaches have been used to estimate potential on a smaller geographical scale. To the best of the authors' knowledge, no rooftop PV potential analysis has been performed using 3D building models with roof geometries for all of Germany.

\subsection{Methodology}\label{sssec:rtpv-met}
In the present study we used 3D building models based on LiDar in the CityGML format with level of detail (LoD)~2 \cite{biljecki_improved_2016} of the Bundesamt für Kartographie und Geodäsie (BKG) \cite{geobasisdaten__geobasis-de__bkg_2021_3d-gebaudemodelle_2021} (high resolution) to estimate rooftop PV potential in Germany (high geographical extent). 
To this end, $93.1$~million roofs were evaluated to estimate Germany's rooftop PV potential.
LoD2 corresponds to simplified building geometries with standardized roof shapes \citep{bezirksregierung_koln_nutzerinformationen_2021}. 
Therefore, the orientation of the roof, i.e., the tilt and azimuth, can be estimated. However, information regarding superstructures reducing the usable area for rooftop PV is missing. 
\par
In order to extract tilt and azimuth from roof geometries, the normal vector on the plane of the geometry is determined. 
The north-based azimuth is defined as the angle between the north vector and x-y part of the normal vector:
\begin{equation}
    azi = \arccos \left( \overrightarrow{n}\cdot 
    [0, 1, 0]^{T} \right)
\end{equation}
The tilt can be determined by calculating the angle between the z-vector and normal vector:
\begin{equation}
    tilt = \arccos \left( \overrightarrow{n}\cdot 
    [0, 0, 1]^{T}
    \right)
\end{equation}
To further estimate the potential to install PV modules on rooftops, in a first step the area of the roof geometries was directly retrieved from the data. 
However, the usable area for rooftop PV is limited by shading or obstacles like windows or chimneys, which is taken into account by incorporating a factor. 
In this paper, the factor was set to $fac_{area}=0.6$, which is in accordance with literature values.
The capacity potential of individual roofs is estimated by placing modules with an efficiency of $\eta=0.22$ (see \autoref{sssec:ofpv-met}). 
Moreover, roofs with tilts lower than \SI{10}{\degree} are assumed to be flat. 
Modules on flat roofs are placed in a southerly direction with an optimal tilt angle in accordance with \citet{ryberg_generation_2019}.
For row spacing on flat roofs, an additional factor of $fac_{RS}=0.5$ is employed. 
The capacity calculation for single roofs can be summarized as follows:
\begin{subnumcases}{P_{PV, peak}=}
    0.6 \cdot \eta_{PV} \cdot A_{roof}, & for tilt $\geq$ \SI{10}{\degree} \\
    0.6 \cdot 0.5 \cdot \eta_{PV} \cdot A_{roof}, & for tilt $<$ \SI{10}{\degree}
\end{subnumcases}
Areas, that correspond to less than \SI{1}{\kilo\watt_p} capacity are excluded from the potential. 
This translates to an area threshold of \SI{7.6}{\metre\squared} for tilted roofs and \SI{15.2}{\metre\squared} for flat ones.
The potentials in \textbf{trep-db} are published in groups for each municipality, nuts3-region, and federal state, but not individually for each building due to the large overhead storage.
To this end, a fixed grouping method is used:
One flat group includes all items up to a tilt angle of 20°. 
Eight~azimuth groups (North (N), North-West (NW), West (W), South-West (SW), South (S), South-East (SE), East (E), and North-East (NE)) contain all elements from 20-90° in the respective directions. 
For rooftop PV, no scenarios are considered, but the \textbf{trep-db} includes two datasets, representing two scenarios:
\begin{itemize}
    \item All roofs  
    \item No northern roofs: Exclusion of north facing groups (N,NW,NE)
\end{itemize}

\subsection{Results \& Discussion}\label{sssec:rtpv-res}
Based on the presented workflow, a potential of \SI{625}{\giga\watt_p} is estimated for all rooftops in Germany.
However, the potential decreases to \SI{492}{\giga\watt_p} if north-facing groups are excluded.  
\par
\begin{figure}[h]
	\centering
		\includegraphics[width=\columnwidth]{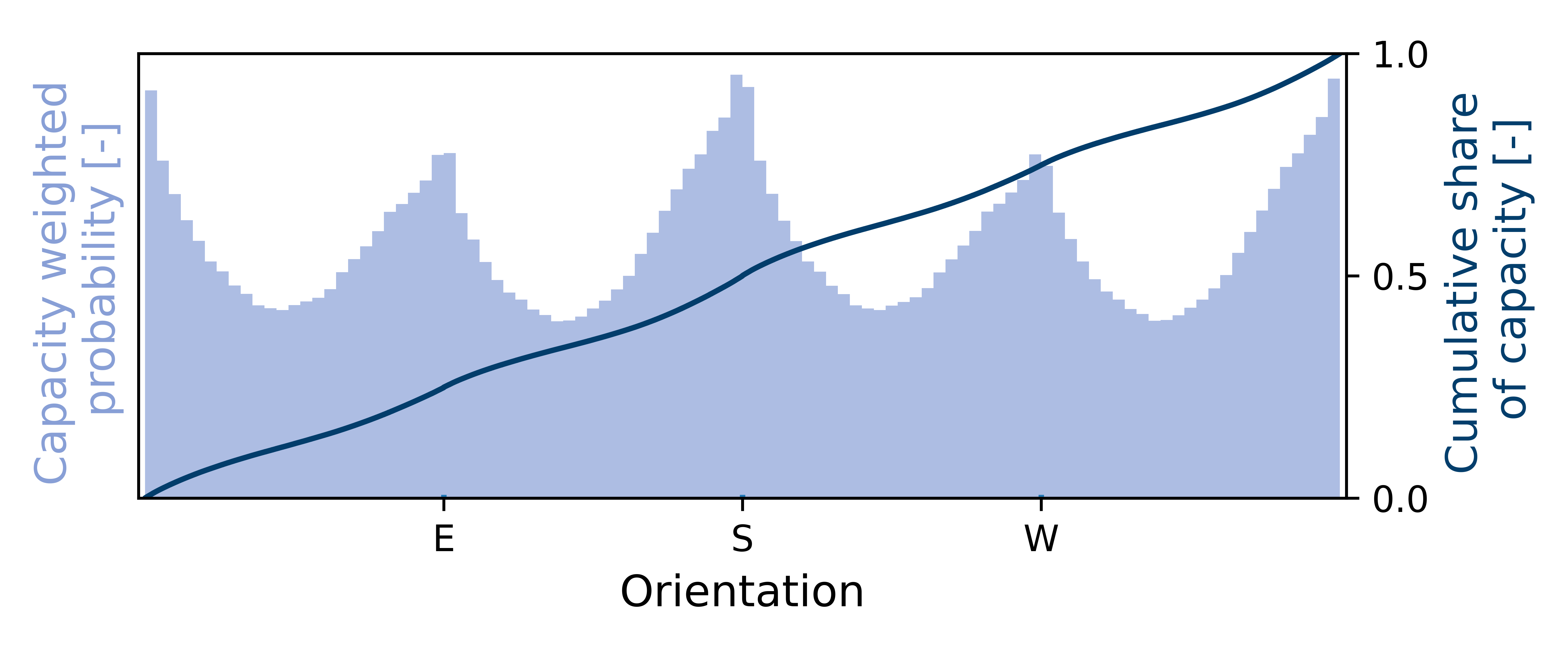}
	  \caption{Capacity weighted azimuth distribution of German roofs.}\label{fig:azimuth_dist}
\end{figure}
\begin{figure}[h]
	\centering
		\includegraphics[width=\columnwidth]{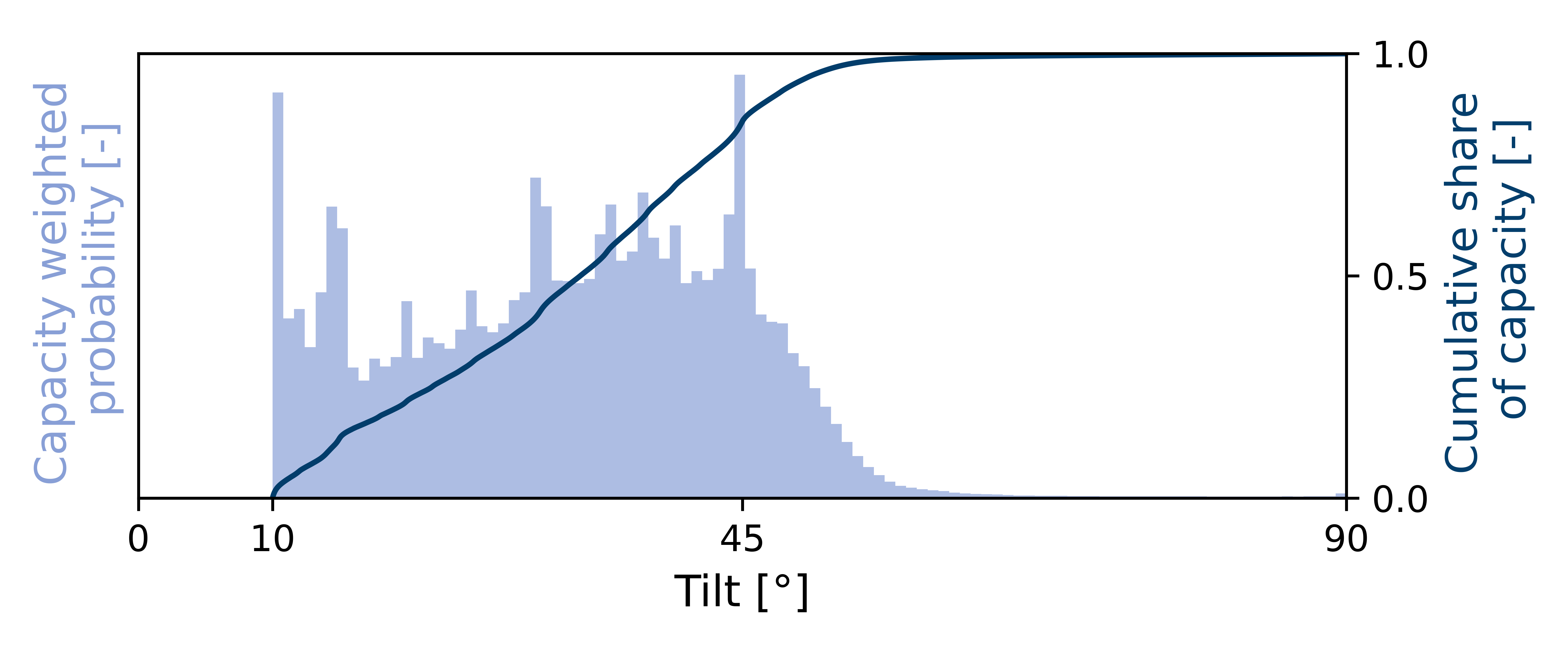}
	  \caption{Capacity-weighted tilt distribution of German roofs.}\label{fig:tilt_dist}
\end{figure}

\autoref{fig:azimuth_dist} and \autoref{fig:tilt_dist} demonstrate the capacity-weighted potential distribution of the azimuth angle and tilt angle for saddle roofs. 
The cardinal directions show the highest capacity potential. 
A substantial capacity potential is located on northern facing roofs.
\autoref{fig:tilt_dist} shows the distribution of the capacity weighted potential distribution of the tilt angle. 
As only saddle roofs are considered for the plots, only tilt angles from 10° upwards are present.
The main capacity potentials can be observed up to a tilt angle of 45° (see \autoref{fig:tilt_dist}).

To validate the presented workflow, the results for the federal state North Rhine-Westphalia are compared against \citet{grothues_landesweite_2018}.
When aligning the workflows by lowering the efficiency to $17\%$, using $f_{RS}=0.4$ and applying the same limit for areas (\SI{7}{\metre\squared} for tilted roofs and \SI{17.5}{\metre\squared} for flat ones) the capacity of roofs in North Rhine-Westphalia adds up to \SI{84.0}{\giga\watt_p}, which is comparable to the \SI{81.4}{\giga\watt_p} estimated by \citet{grothues_landesweite_2018}. 
It should be noted that \citet{grothues_landesweite_2018} exclude rooftops with an irradiation lower than \SI{814}{{\kilo\watt\hour\per\metre\squared}}, which is not reproduced, and therefore a higher result is to be expected. 
Furthermore, the potential reduction due to roof obstructions is not comparable in the two studies. 

\par
\begin{figure}[h]
	\centering
\includegraphics[width=0.50\textwidth]{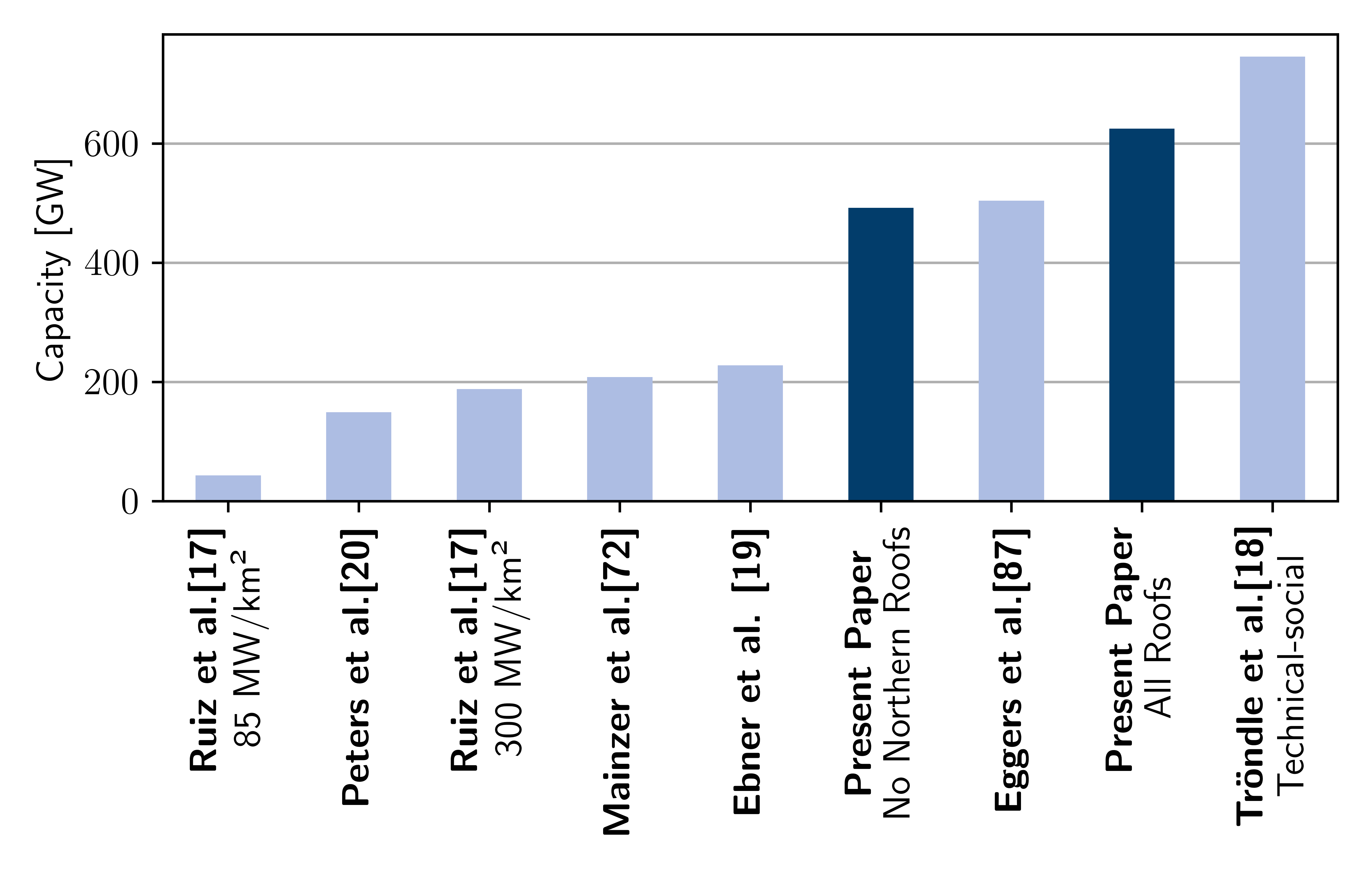}
	  \caption{Comparison of rooftop PV potentials for studies providing capacity potential at the national level.}\label{fig:comparison_rtpv}
\end{figure}

\autoref{fig:comparison_rtpv} presents the Germany-wide capacity potential of rooftop PV in this paper and other studies. 
The comparison shows a range of \SI{43}{\giga\watt_p}~\citep{ruiz_enspreso_2019} to \SI{746}{\giga\watt_p}~\citep{trondle_home-made_2019}.
However, comparability between the studies is limited.
The studies significantly differ in their methodologies for rooftop area estimations and assumptions as the reduction factor for superstructures and efficiency in the capacity estimation.
Furthermore, some studies only consider the rooftop PV potential for individual building categories, e.g., residential in \citet{mainzer_assessment_2017}. The only capacity potential, which has been estimated by a high resolution approach~(see \autoref{sssec:rtpv-lit}) assumes \SI{504}{\giga\watt}~\citep{eggers_pv-ausbauerfordernisse_2020}, which is of similar magnitude to the potentials estimated in our approach.

\par
In future, the estimation of unusable areas on roofs can be improved when LoD3/LoD4-data becomes available for Germany. Alternatively, as proposed by \citet{walch_big_2020}, the LoD4-Dataset from Geneva could be used to estimate the share of unusable roof areas by means of a machine learning approach.

\section{Discussion}\label{sec:disc}
\begin{figure}[h]
	\centering
\includegraphics[width=0.50\textwidth]{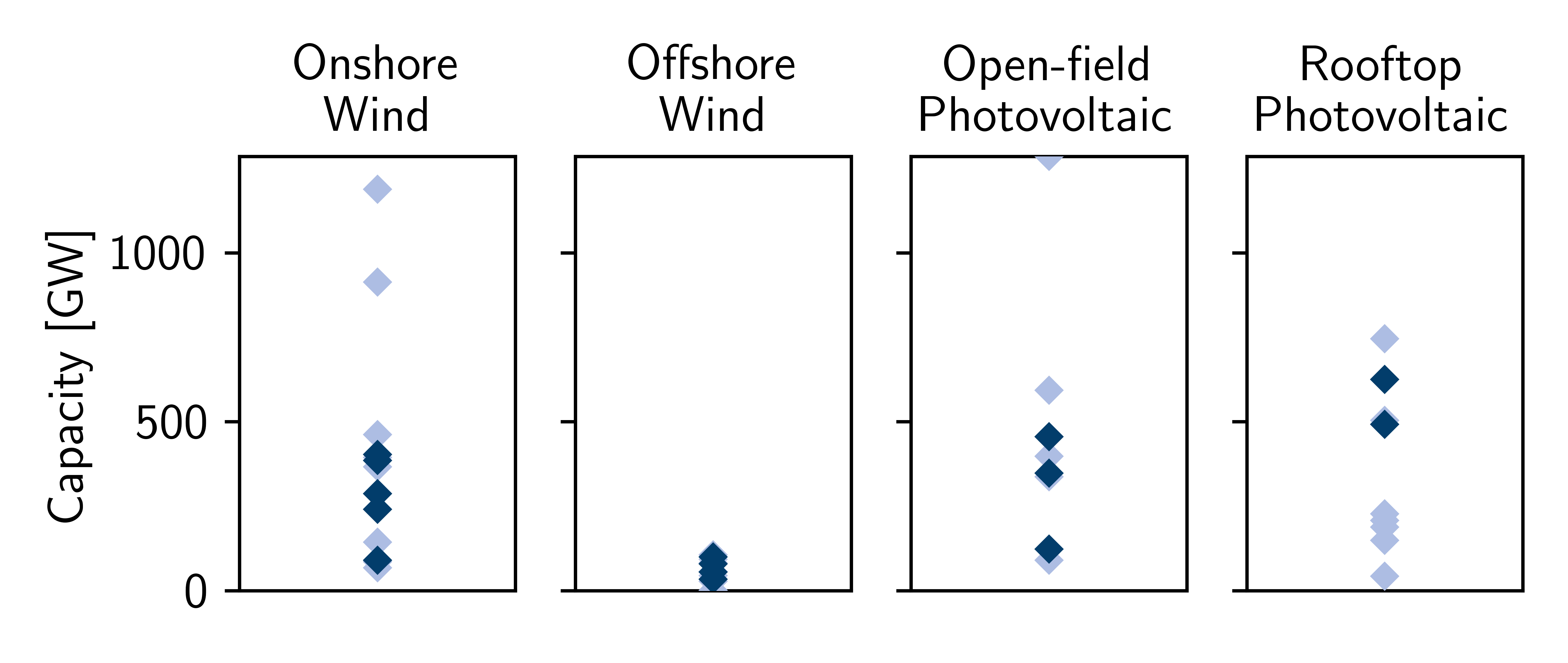}
	  \caption{Potential capacities of different potential studies on the national level for onshore wind, offshore wind, open-field photovoltaic, and rooftop photovoltaic.}\label{fig:violin_plot}
\end{figure}

The capacity estimations of the reviewed literature for different technologies in Germany are shown in \autoref{fig:violin_plot}.
The potentials exhibit large variations between the studies and scenarios. Our data analysis (\autoref{sssec:data}) shows one reason for the present deviations: the used input datasets make the findings of different studies using different datasets hardly comparable. Another reason for the significantly deviating potentials is the chosen exclusion criteria and the corresponding buffers. Neglecting residential buildings in wind potential analyses, for example, can lead to considerable over-estimations. In future work, a sensitivity analysis regarding the criteria could help quantify the effects in a more in-depth manner. Nevertheless, without applying high-quality datasets, the conclusions would not be meaningful. In other regions of the world, OSM may not be as comprehensive and a dataset with official characteristic such as Basis-DLM may not be available. In future work, world-wide high-resolution datasets should be used to reevaluate global renewable potentials. 
\par
\begin{figure*}
	\centering
\includegraphics[width=496pt]{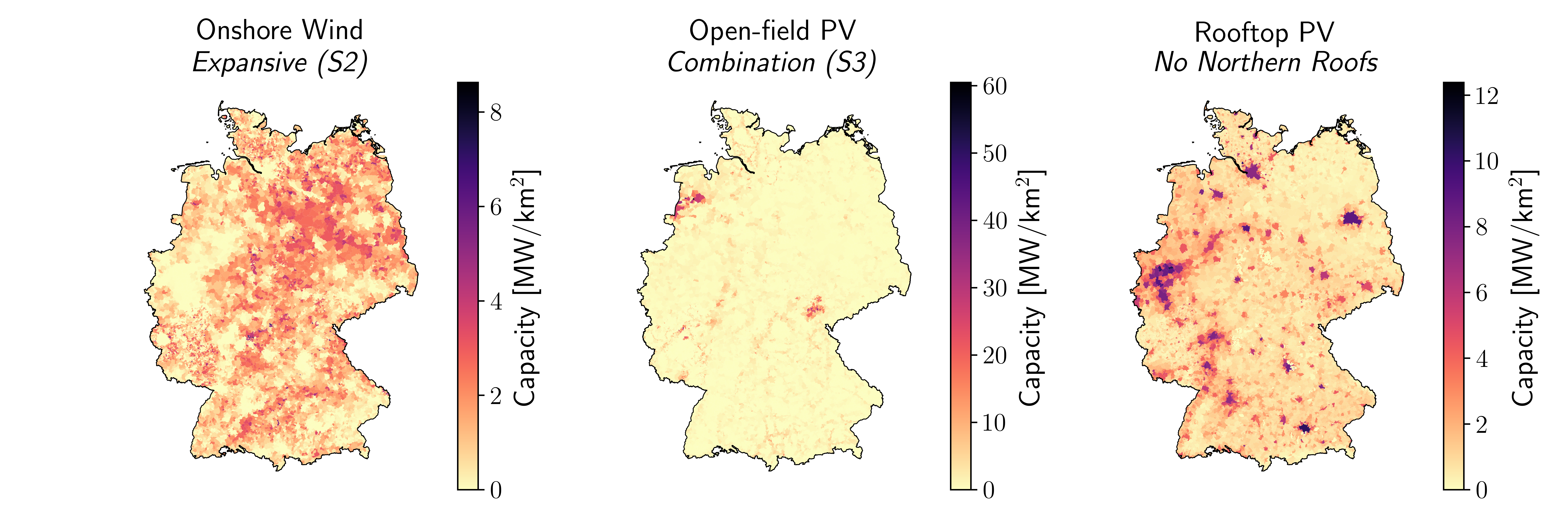}
	  \caption{Potential capacity density in Germany's Municipalities for different technologies.}\label{fig:scenario_maps}
\end{figure*}
\begin{figure*}
	\centering
\includegraphics[width=420pt]{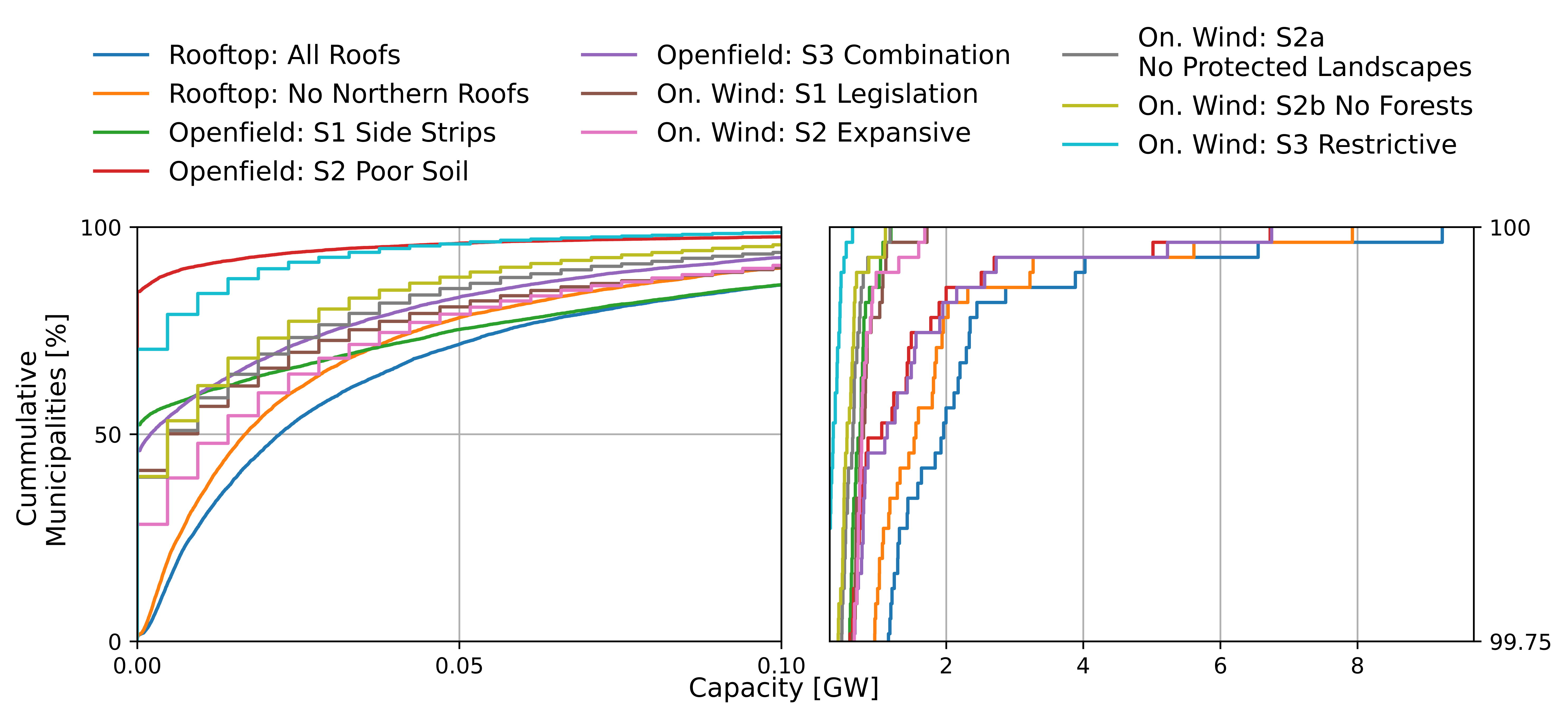}
	  \caption{Cumulative municipalities over capacity potential for the scenarios of open-field, rooftop PV and wind onshore.}\label{fig:share_muni}
\end{figure*}

\autoref{fig:scenario_maps} shows the spatial distribution of the capacity potential in Germany's municipalities for one scenario of onshore wind, open-field PV, and rooftop PV, which visualizes the spatial heterogeneity of the results. 
Additionally, \autoref{fig:share_muni} shows the cumulative share of Germany's municipalities over all scenarios for the regarded technologies and 
the impact of boundary conditions, e.g., pre-selected areas or exclusions, on regional potentials in municipalities.
For the wind scenarios, the share of municipalities without potential highly differs between $28\%$ and $70\%$, showing the high impact of the exclusion criteria on the regional level.
For all open-field PV scenarios, a large number of municipalities are without potentials, whereas rooftop PV potentials are more equally distributed.
As is shown on the right-hand side of the figure, the maximum capacities found in the municipalities highly vary for the technologies and scenarios between \SI{634}{\mega\watt} (On. Wind, S3 Restrictive) and \SI{9.3}{\giga\watt} (Rooftop PV, All Roofs).
For the use case of energy system models, the unequal distribution of potentials in municipalities emphasizes the importance of detailed and qualitative potentials for regional analyses.

\section{Conclusions}\label{sec:conclusion}
The first time evaluation of land use datasets and the proposed renewable potential scenarios reveal significant biases in the datasets that are commonly applied by the energy systems community:
The use of Corine Land Cover \citep{copernicus_programme_corine_2018} leads to a significant overestimation of the potentially usable area for renewable energy technologies by a factor of 4.6 to 5.2 in comparison to Basis-DLM \citep{geobasisdaten__geobasis-de__bkg_2021_digitales_2021} and Open Street Map~\citep{openstreetmap_contributors_open_2021}, respectively. High-quality datasets are needed to supply reliable input for policy-makers or energy system models. In the case of Germany, Basis-DLM \citep{geobasisdaten__geobasis-de__bkg_2021_digitales_2021} and Open Street Map~\citep{openstreetmap_contributors_open_2021} provide similar information for several categories, especially line-like features, such as power lines or railways. For others, Basis-DLM and Open Street Map show significant differences, e.g., industry/commercial, lakes, or rivers.
Furthermore, the impact of several exclusion criteria is apparent in the presented scenarios and variations. For example, the disregard of residential buildings in onshore wind analyses leads to an overestimation of the capacity by up to $34\%$ for our scenarios. \par
The presented scenarios for onshore wind, offshore wind, open-field, and rooftop photovoltaic are published in the open access database \textbf{trep-db}. For rooftop photovoltaic, the potentials are estimated by CityGML Level of Detail 2 data for the first time, which leads to a \SI{492}{\giga\watt_p} potential capacity when northerly-facing roofs are excluded. 
\par
This work should motivate improvement in the internationally-available land use datasets, as their quality has a significant impact on the correct design of regulations for renewable resource emplacements.

\section{Acknowledgements}
This work was supported by the BMWi (German Federal Ministry of Economic
Affairs and Energy) [promotional reference 3EE5031D]; and the Helmholtz Association under the program “Energy System Design.”


\newpage


\bibliography{PotentialsPaper}



\end{document}